
\def\LaTeX{%
  \let\Begin\begin
  \let\End\end
  \let\salta\relax
  \let\finqui\relax
  \let\futuro\relax}

\def\UK{\def\our{our}\let\sz s}
\def\USA{\def\our{or}\let\sz z}

\UK



\LaTeX

\USA


\salta

\documentclass[twoside,12pt]{article}
\setlength{\textheight}{24cm}
\setlength{\textwidth}{16cm}
\setlength{\oddsidemargin}{2mm}
\setlength{\evensidemargin}{2mm}
\setlength{\topmargin}{-15mm}
\parskip2mm


\usepackage[usenames,dvipsnames]{color}
\usepackage{amsmath}
\usepackage{amsthm}
\usepackage{amssymb}
\usepackage[mathcal]{euscript}
\usepackage{cite}
%
%
%


\definecolor{viola}{rgb}{0.3,0,0.7}
\definecolor{ciclamino}{rgb}{0.5,0,0.5}
\definecolor{rosso}{rgb}{0.85,0,0}

\def\pier #1{#1}
\def\juerg #1{#1}
\def\gian #1{#1}
\def\gianni #1{#1}



\bibliographystyle{plain}


%

\finqui

\def\Beq{\Begin{equation}}
\def\Eeq{\End{equation}}
\def\Bsist{\Begin{eqnarray}}
\def\Esist{\End{eqnarray}}

\def\Bthm{\Begin{theorem}}
\def\Ethm{\End{theorem}}
\def\Blem{\Begin{lemma}}
\def\Elem{\End{lemma}}

\def\Bcor{\Begin{corollary}}
\def\Ecor{\End{corollary}}
\def\Brem{\Begin{remark}\rm}
\def\Erem{\End{remark}}

\def\Bcenter{\Begin{center}}
\def\Ecenter{\End{center}}
\let\non\nonumber
\def\proof{\noindent{\sc Proof}: \quad}




\def\step #1 \par{\medskip\noindent{\bf #1.}\quad}
\def\jstep #1: \par {\vspace{2mm}\noindent\underline{\sc #1 :}\par\nobreak\vspace{1mm}\noindent}


\def\Holder{H\"older}
\def\Frechet{Fr\'echet}
\def\aand{\quad\hbox{and}\quad}

\def\lhs{left-hand side}
\def\rhs{right-hand side}
\def\sfw{straightforward}



\def\multibold #1{\def\arg{#1}%
  \ifx\arg\pto \let\next\relax
  \else
  \def\next{\expandafter
    \def\csname #1#1#1\endcsname{{\bf #1}}%
    \multibold}%
  \fi \next}

\def\pto{.}

\def\multical #1{\def\arg{#1}%
  \ifx\arg\pto \let\next\relax
  \else
  \def\next{\expandafter
    \def\csname cal#1\endcsname{{\cal #1}}%
    \multical}%
  \fi \next}


\def\multimathop #1 {\def\arg{#1}%
  \ifx\arg\pto \let\next\relax
  \else
  \def\next{\expandafter
    \def\csname #1\endcsname{\mathop{\rm #1}\nolimits}%
    \multimathop}%
  \fi \next}

\multibold
qwertyuiopasdfghjklzxcvbnmQWERTYUIOPASDFGHJKLZXCVBNM.

\multical
QWERTYUIOPASDFGHJKLZXCVBNM.

\multimathop
diag dist div dom mean meas sign supp .

\def\Span{\mathop{\rm span}\nolimits}


\def\accorpa #1#2{\eqref{#1}--\eqref{#2}}
\def\Accorpa #1#2 #3 {\gdef #1{\eqref{#2}--\eqref{#3}}%
  \wlog{}\wlog{\string #1 -> #2 - #3}\wlog{}}


\def\separa{\noalign{\allowbreak}}

\def\somma #1#2#3{\sum_{#1=#2}^{#3}}

\def\graffe #1{\mathopen\{#1\mathclose\}}

\def\<#1>{\mathopen\langle #1\mathclose\rangle}
\def\norma #1{\mathopen \| #1\mathclose \|}

\def\iot {\int_0^t}
\def\ioT {\int_0^T}
\def\inttT {\int_t^T}
\def\intsT {\int_s^T}
\def\intQt{\int_{Q_t}}
\def\intQ{\int_Q}
\def\iO{\int_\Omega}
\def\iG{\int_\Gamma}
\def\intS{\int_\Sigma}
\def\intSt{\int_{\Sigma_t}}

\def\bintQt{\int_{Q^t}}
\def\bintSt{\int_{\Sigma^t}}
\def\bintQs{\int_{Q^s}}
\def\bintSs{\int_{\Sigma^s}}

\def\dt{\partial_t}
\def\dn{\partial_\nu}

\def\cpto{\,\cdot\,}

\def\checkmmode #1{\relax\ifmmode\hbox{#1}\else{#1}\fi}
\def\aeO{\checkmmode{a.e.\ in~$\Omega$}}
\def\aeQ{\checkmmode{a.e.\ in~$Q$}}
\def\aeG{\checkmmode{a.e.\ on~$\Gamma$}}
\def\aeS{\checkmmode{a.e.\ on~$\Sigma$}}
\def\aet{\checkmmode{a.e.\ in~$(0,T)$}}

\def\aat{\checkmmode{for a.a.~$t\in(0,T)$}}


\def\erre{{\mathbb{R}}}
\def\erren{\erre^n}

\def\enne{{\mathbb{N}}}

\def\VV{{\mathbb{V}}}
\def\HH{{\mathbb{H}}}
\def\AA{{\mathbb{A}}}
\def\BB{{\mathbb{B}}}
\def\CC{{\mathbb{C}}}
\def\II{{\mathbb{I}}}




\def\genspazio #1#2#3#4#5{#1^{#2}(#5,#4;#3)}
\def\spazio #1#2#3{\genspazio {#1}{#2}{#3}T0}

\def\L {\spazio L}
\def\H {\spazio H}
\def\W {\spazio W}

\def\C #1#2{C^{#1}([0,T];#2)}



\def\Lx #1{L^{#1}(\Omega)}
\def\Hx #1{H^{#1}(\Omega)}

\def\LxG #1{L^{#1}(\Gamma)}
\def\HxG #1{H^{#1}(\Gamma)}

\def\LQ #1{L^{#1}(Q)}
\def\LS #1{L^{#1}(\Sigma)}

\def\Ldue{\Lx 2}

\def\Huno{\Hx 1}
\def\Hdue{\Hx 2}

\def\HunoG{\HxG 1}
\def\HdueG{\HxG 2}

\def\LdueG{\LxG 2}


\def\LQ #1{L^{#1}(Q)}


\let\theta\vartheta
\let\eps\varepsilon
\let\phi\varphi

\let\TeXchi\chi                         
\newbox\chibox
\setbox0 \hbox{\mathsurround0pt $\TeXchi$}
\setbox\chibox \hbox{\raise\dp0 \box 0 }
\def\chi{\copy\chibox}


\def\QED{\hfill $\square$}


\def\CO{C_\Omega}

\def\suG{_{|\Gamma}}

\def\VG{V_\Gamma}
\def\HG{H_\Gamma}
\def\WG{W_\Gamma}
\def\nablaG{\nabla_\Gamma}
\def\DeltaG{\Delta_\Gamma}
\def\muG{\mu_\Gamma}
\def\rhoG{\rho_\Gamma}
\def\tauO{\tau_\Omega}
\def\tauG{\tau_\Gamma}
\def\fG{f_\Gamma}
\def\vG{v_\Gamma}
\def\wG{w_\Gamma}
\def\zG{z_\Gamma}

\def\xiG{\xi_\Gamma}
\def\etaG{\eta_\Gamma}

\def\Mu{(\mu,\muG)}
\def\Rho{(\rho,\rhoG)}

\def\Xi{(\xi,\xiG)}
\def\psiG{\psi_\Gamma}
\def\pG{p_\Gamma}
\def\qG{q_\Gamma}

\def\QG{Q_\Gamma}

\def\rhoz{\rho_0}

\def\ueps{u^\eps}

\def\peps{p^\eps}
\def\pGeps{p_\Gamma^\eps}
\def\qeps{q^\eps}
\def\qGeps{q_\Gamma^\eps}
\def\phieps{\phi^\eps}

\def\Pi{\hat\pi}

\def\etan{\eta^n}
\def\xin{\xi^n}
\def\etaGn{\eta_\Gamma^n}
\def\xiGn{\xi_\Gamma^n}
\def\Etan{(\etan,\etaGn)}
\def\Xin{(\xin,\xiGn)}

\def\un{u_n}
\def\mun{\mu_n}
\def\muGn{\mu_{n_\Gamma}}
\def\rhon{\rho_n}
\def\rhoGn{\rho_{n_\Gamma}}

\def\ei{e^i}
\def\eGi{e_\Gamma^i}
\def\Ei{(\ei,\eGi)}
\def\ej{e^j}
\def\eGj{e_\Gamma^j}
\def\Ej{(\ej,\eGj)}

\def\rhomin{\rho_*}
\def\rhomax{\rho^*}

\def\calVp{\calV^{\,*}}
\def\calVn{\calV_n}

\def\normaHH #1{\norma{#1}_{\calH}}
\def\normaVV #1{\norma{#1}_{\calV}}

\def\ei{e^i}
\def\ej{e^j}
\def\eGi{e_\Gamma^i}
\def\eGj{e_\Gamma^j}

\let\hat\widehat

\def\Uad{\calU_{ad}}

\def\mub{\overline\mu}
\def\muGb{\overline\mu_\Gamma}
\def\rhob{\overline\rho}
\def\rhoGb{\overline\rho_\Gamma}
\def\ub{\overline u}

\def\hmu{\hat\mu}
\def\hrho{\hat\rho}
\def\hmuQ{\hmu_Q}
\def\hmuS{\hmu_\Sigma}
\def\hrhoQ{\hrho_Q}
\def\hrhoS{\hrho_\Sigma}
\def\hrhoO{\hrho_\Omega}
\def\hrhoG{\hrho_\Gamma}

\def\VVp{{\VV\,}^*}

\Begin{document}


%
\title{Optimal velocity control of a viscous\\[2mm]
Cahn--Hilliard system with convection\\[2mm] 
and dynamic boundary conditions\footnote{This work received 
a partial support from the MIUR-PRIN Grant 2015PA5MP7 
``Calculus of Variations'', the GNAMPA (Gruppo Nazionale per l'Analisi
Matematica, la Probabilit\`{a} e loro Applicazioni) of INDAM (Istituto Nazionale
di Alta Matematica) and the IMATI -- C.N.R. Pavia for PC and GG.}}
\author{}
\date{}
\maketitle
\Bcenter
\vskip-2cm
{\large\sc Pierluigi Colli$^{(1)}$}\\
{\normalsize e-mail: {\tt pierluigi.colli@unipv.it}}\\[.25cm]
{\large\sc Gianni Gilardi$^{(1)}$}\\
{\normalsize e-mail: {\tt gianni.gilardi@unipv.it}}\\[.25cm]
{\large\sc J\"urgen Sprekels$^{(2)}$}\\
{\normalsize e-mail: {\tt sprekels@wias-berlin.de}}\\[.45cm]
$^{(1)}$
{\small Dipartimento di Matematica ``F. Casorati'', Universit\`a di Pavia}\\
{\small and Research Associate at the IMATI -- C.N.R. Pavia}\\
{\small via Ferrata 5, 27100 Pavia, Italy}\\[.2cm]
$^{(2)}$
{\small Department of Mathematics}\\
{\small Humboldt-Universit\"at zu Berlin}\\
{\small Unter den Linden 6, 10099 Berlin, Germany}\\[2mm]
{\small and}\\[2mm]
{\small Weierstrass Institute for Applied Analysis and Stochastics}\\
{\small Mohrenstrasse 39, 10117 Berlin, Germany}
\Ecenter
\Begin{abstract}
\noindent In this paper, we investigate a distributed optimal control problem for a 
convective viscous Cahn--Hilliard system with dynamic boundary conditions.
Such systems govern phase separation processes between two phases taking place in an incompressible
fluid in a container and, at the same time, on the container boundary.  
The cost functional is of standard tracking type, while the control is exerted
by the velocity of the fluid in the bulk. 
In this way, the coupling between the state (given by the associated order parameter and chemical potential) and 
control variables in the governing system of nonlinear partial differential equations is bilinear, 
which presents an additional difficulty
for the analysis. The nonlinearities in the bulk and surface free energies are of 
logarithmic type, which entails that the thermodynamic forces driving the phase separation process
may become singular. We show existence for the optimal control problem under investigation,
prove the Fr\'echet differentiability of the associated control-to-state mapping in suitable Banach
spaces, and derive the first-order necessary optimality conditions in terms of a variational
inequality and the associated adjoint system. Due to the strong nonlinear couplings between
state variables and control, the corresponding proofs require a considerable analytical effort.
     

\End{abstract}
\salta
\pagestyle{myheadings}
\newcommand\testopari{\sc Colli \ --- \ Gilardi \ --- \ Sprekels}
\newcommand\testodispari{\sc Velocity control of convective Cahn--Hilliard system}
\markboth{\testopari}{\testodispari}
\finqui
%

\section{Introduction}
\label{INTRO}
\setcounter{equation}{0}

Let $\Omega\subset\erre^3$ denote some open, bounded and connected set having a smooth boundary $\Gamma$
and unit outward normal $\,\nu$.
We denote by $\dn$, $\nablaG$, $\Delta_\Gamma$ the outward normal derivative, the tangential gradient,
and the Laplace--Beltrami operator on $\Gamma$, in this order. Moreover, we fix some final time $T>0$ and
introduce for every $t\in (0,T]$ the sets $Q_t:=\Omega\times (0,t)$ and $\Sigma_t:=\Gamma\times (0,t)$,
where we put, for the sake of brevity, $Q:=Q_T$ and $\Sigma:=\Sigma_T$. We then consider the following optimal control
problem:

\vspace{3mm}\noindent
{\bf (CP)} \quad Minimize the cost functional
\begin{align}
\label{cost} 
  & \calJ(\mu,\muG,\rho,\rhoG,u)
  := \frac{\beta_1}2 \intQ |\mu-\hmuQ|^2
  + \frac{\beta_2}2 \intS |\muG-\hmuS|^2
  \non
  \\
  & \quad {}
  + \frac{\beta_3}2 \intQ |\rho-\hrhoQ|^2
  + \frac{\beta_4}2 \intS |\rho-\hrhoS|^2
  \non
  \\
  & \quad {}
  + \frac{\beta_5}2 \iO |\rho(T)-\hrhoO|^2
  + \frac{\beta_6}2 \iG |\rhoG(T)-\hrhoG|^2
  + \frac{\beta_7}2 \intQ |u|^2\,,
\end{align}
subject to  the state system  
\begin{align}
\label{ss1}
&\dt\rho+\nabla\rho\cdot u-\Delta\mu=0 \quad\mbox{in }\,Q\,,\\[1mm]
\label{ss2}
&\tauO\dt\rho-\Delta\rho+f'(\rho)=\mu \quad\mbox{in }\,Q\,,\\[1mm]
\label{ss3}
&\dt\rhoG+\dn\mu-\DeltaG\muG=0 \quad\mbox{and}\quad \mu_{|\Sigma}=\muG\quad\mbox{on }\,\Sigma\,,\\[1mm]
\label{ss4}
&\tauG\dt\rhoG+\dn\rho-\DeltaG\rhoG+\fG'(\rhoG)=\muG \quad\mbox{and}\quad \rho_{|\Sigma}=\rhoG
\quad\mbox{on }\,\Sigma\,,\\[1mm]
\label{ss5}
&\rho(0)=\rho_0\quad\mbox{in }\,\Omega,\quad \rhoG(0)=\rho_{0|\Gamma}\quad\mbox{on }\,\Gamma\,, 
\end{align}
and to the control constraint
\begin{equation}
\label{concon}
u\in\Uad\,,
\end{equation}
\gianni{where $\Uad$ is a suitable closed, convex, and bounded subset of the control space $\calX$ defined~by
\Bsist
  && \calX := \L2Z \cap (L^\infty(Q))^3 \cap (\H1{\Lx3})^3 \pier{,}
  \label{defcalX}
  \\
  \noalign{\noindent where}
  && Z := \graffe{ w\in(\Lx2)^3:\ \div w=0 \ \hbox{in $\Omega$ \ and }\ w\cdot\nu=0 \ \hbox{on $\Gamma$}}.
  \label{defZ}
\Esist
}%
\gianni{In \eqref{cost}}, the constants $\beta_i$, $1\le i\le 7$, are nonnegative but not all zero, and $\widehat\mu_Q$,
$\widehat\mu_\Sigma$, $\widehat\rho_Q$, $\widehat \rho_\Sigma$, $\widehat\rho_\Omega$, and
$\widehat\rho_\Gamma$, are given target functions\gianni.
 We note that the state system \eqref{ss1}--\eqref{ss5} can be seen as a phase field model for a
 phase separation process taking place in an incompressible fluid in the container $\Omega$ and
 on the container boundary~$\Gamma$. In this connection, the variables $\Mu$ and $\Rho$ stand
 for the chemical potential and the order parameter (usually the density of one of the involved
 phases, normalized in such a way as to attain its values in the interval $[-1,1]$) 
 of the phase separation process in the bulk and on the surface, respectively. It is worth
 noting that the total mass of the order parameter is conserved during the separation process;
 indeed, integrating \eqref{ss1} for fixed $t\in (0,T]$ over $\Omega$, and using
\gianni{the condition $u(t)\in Z$} and~\eqref{ss3}, we readily find that
 \begin{equation}
 \label{conserve}
 \dt\Big(\iO\rho(t)+\iG\rhoG(t)\Big)=0\,.
 \end{equation} 
 We also assume that the densities of the local free bulk energy $f$ and the
local free surface energy $\fG$ are of logarithmic type, where the latter dominates the former
in a sense to be made precise later. In the simplest case, we have 
\begin{equation}
\label{flog}
f(r)\pier{{}\simeq{}}\fG(r) \pier{{}\simeq{}} \widehat c_1((1+r)\,\ln(1+r)+(1-r)\ln(1-r))-\widehat c_2r^2, \quad \pier{r\in (-1,1)},
\end{equation}   
\pier{with constants (not necessarily the same) $\widehat c_1>0$ and $\widehat c_2>0$ such that both} $f$ and $\fG$ are nonconvex. Notice 
that the derivatives $f'$ and $\fG'$ are singular at the endpoints $r=\pm 1$.

While there are numerous \pier{contributions} (which cannot be cited here) in the literature that address the questions of well-posedness and asymptotic behavior
for various types (viscous or nonviscous, local or nonlocal, zero Neumann boundary conditions or
dynamic boundary conditions) of Cahn--Hilliard systems, there are
still but a few papers dealing with the associated optimal control problems. In this connection, we
refer to  \cite{CGRS1,CGS2,CGS4,hw,WaNa,ZW} for the case of Dirichlet or zero 
Neumann boundary conditions and
to \pier{\cite{CFGS1,CFGS2,CGS1,CGS3,CGS3bis,FY}} for the case of dynamic boundary conditions. 

Recently, a  rigorous analysis  
for convective Cahn--Hilliard systems has been given 
in \cite{ZL1} for the one-dimensional and in \cite{ZL2} for the two-dimensional case. In \cite{FRS},
the distributed optimal control of a two-dimensional Cahn--Hilliard/Navier--Stokes system was 
analyzed. We also mention the papers \cite{HHKW,HW3,HW1,HW2}, which deal with the optimal control of three-dimensional 
Cahn--Hilliard/Navier--Stokes systems, however in the time-discretized version.

A distinguishing feature of this paper is that we use the fluid velocity as the control variable in the
convective Cahn--Hilliard system. In practice, this can be realized by placing either a mechanical stirring device or 
an ultrasound emitter into the container. Another option is, in the case of electrically conducting fluids
like molten metals, to make use of magnetic fields (for such an application, see \cite{Kudla}). 
To the authors' best knowledge, the only existing mathematical contribution, in which the fluid velocity
is used as the control in a convective Cahn--Hilliard system in three dimensions of space, is the recent contribution
\cite{RS}. In comparison with the situation investigated in \cite{RS}, the main novelties of our paper  
are the following: while in \cite{RS} a nonlocal convective Cahn--Hilliard system with a possibly
degenerating mobility and zero
Neumann boundary conditions was studied, we consider here a viscous local Cahn--Hilliard system with constant 
mobility (normalized to unity) and the more difficult dynamic
boundary conditions. In the recent paper \cite{CGS13}, rather general and strong well-posedness results
for this situtation have been established \pier{(we also like to quote the contributions \cite{CF1, CF2} for the nonconvective case).} 

In our analysis, we will take advantage of the results shown in \cite{CGS13}. It turns out that the bilinear
coupling between control and state makes it necessary to allow only controls $u$ which, among other constraints,
have to obey the somewhat unusual regularity condition $u\in H^1(0,T;L^3(\Omega)^3)$. But, as a matter of
fact, this is exactly the kind of regularity that guarantees the existence of a unique solution to the state
system having sufficient regularity properties. 
Under these premises, we will be able to show the Fr\'echet differentiability 
of the control-to-state operator in suitable Banach spaces. 
\gian{Finally, we can prove the existence of an optimal control
and, in a slightly less general setting,
we also derive proper first-order necessary conditions for optimality.}

The paper is organized as follows: in the following Section~\ref{STATE}, 
we state the general assumptions for our problem, 
and we collect known results for the state system \eqref{ss1}--\eqref{ss5}. 
Section~\ref{FRECHET} brings an analysis of the differentiability properties of the control-to-state mapping, 
while in Section~\ref{ADJOINT} 
\gian{we prove existence and the first-order necessary optimality conditions for the control problem.}

Throughout this paper, we will denote for a general Banach space $X$ by $\|\cdot\|_X$ its norm and by
$X^*$ its dual space. Moreover, $\langle\cdot,\cdot\rangle_X$ denotes the dual pairing between $X^*$ and $X$.
The only exception from this convention for the norms is given by the spaces $L^p$ constructed on
$\Omega$, $\Gamma$, $Q$, $\Sigma$ and their powers, for $1\le p\le\infty$, whose norms will be denoted by
$\|\cdot\|_p$. We will also  repeatedly use Young's inequality
\begin{equation}
\label{Young}
a\,b\le \delta\,|a|^2+\mbox{$\frac 1{4\delta}$}\,|b|^2 \quad\mbox{for all } \,a,b\in\erre\quad\mbox{and}\quad
\delta>0,
\end{equation}
as well as the continuity of the embeddings $H^1(\Omega)\subset L^p(\Omega)$ for $1\le p\le 6$ and 
$\Hdue\subset C^0(\overline\Omega)$. Notice that the latter embedding is also compact, while this holds true
for the former embeddings only if $p<6$.


\section{General assumptions and the state system}
\label{STATE}
\setcounter{equation}{0}

In this section, we introduce the general setting of our control 
problem and state some known results on the state system 
\eqref{ss1}--\eqref{ss5}. To begin with, we introduce the spaces
\begin{align}
  & H := \Ldue \,, \quad  
  V := \Huno 
  \aand
  W := \Hdue,
  \label{defspaziO}
  \\
  & \HG := \LdueG \,, \quad 
  \VG := \HunoG 
  \aand
  \WG := \HdueG,
  \label{defspaziG}
  \\
  & \calH := H \times \HG \,, \quad
  \calV := \graffe{(v,\vG) \in V \times \VG : \ \vG = v\suG}
  \non
  \\
  & \aand
  \calW := \bigl( W \times \WG \bigr) \cap \calV \,.
  \label{defspaziprod}
\end{align}
\gian{Moreover, we recall the definition~\eqref{defcalX} of~$\calX$.}

We make the following assumptions on the data of our problem:
\begin{description}
\item{(A1)} \quad$(\rho_0,\rho_{0|\Gamma}) \in\calW$, and we have $\,-1<\rho_0(x)<1\,$ for all $\,x\in\overline\Omega$. 
\item{(A2)} \quad $\tauO>0$ and $\tauG>0$.
\item{(A3)} \quad $f,\fG\in C^3(-1,1)$ can be written as $\,f=f_1+f_2$ and $\fG=f_{\Gamma 1}+f_{\Gamma 2}$, 
 where\linebreak
\hspace*{5mm}$\,f_2,f_{\Gamma 2}\in C^3[-1,1]$ and  
\begin{align}
\label{domina}
&\exists \ \gamma_1>0,\,\gamma_2>0: \quad|f_1'(r)|\le \gamma_1|f_{\Gamma 1}'(r)|+\gamma_2 \quad\forall\,r\in (-1,1),\\[1mm]
\label{singul}
&\lim_{r\searrow -1} f_1'(r)=\lim_{r\searrow -1} f_{\Gamma 1}'(r)=-\infty, \qquad
\lim_{r\nearrow +1} f_1'(r)=\lim_{r\nearrow +1} f_{\Gamma 1}'(r)=+\infty. 
\end{align}
\item{(A4)} \quad The constants $\,\beta_i$, $1\le i\le 7$, are all nonnegative but not all
 equal to zero, and \hspace*{5mm} it holds $\,\widehat \rho_Q, \widehat\mu_Q\in L^2(Q)$, $\,\widehat \rho_\Sigma,
 \widehat\mu_\Sigma\in L^2(\Sigma)$, $\,\widehat\rho_\Omega\in\Ldue$, and $\,\widehat\rho_\Gamma\in L^2(\Gamma)$. 
\gian{%
\item{(A5)} \quad
The function $\overline U\in L^\infty(Q)$ and the constant $R_0>0$ make the 
admissible set
\begin{align}
  & \Uad := \gianni{\bigl\{ u\in\calX: \ |u|\leq\overline U \ \aeQ, \ \norma u_{\calX}\leq R_0 \bigr\}}
  \label{defUad}
\end{align}
\hspace*{5mm}nonempty.}
\end{description}

\noindent For the following analysis, it is convenient to fix once and for all some open ball in $\calX$ that contains $\Uad$. We therefore assume:
\begin{description}
\item{(A6)} \quad Let $R>0$ be fixed such that $\,\,\Uad\subset{\cal U}_R:=\{u\in\gian{\calX}: \,\|u\|_\calX< R\}$.
\end{description} 

\Brem
The condition \eqref{domina} means, loosely speaking, that the thermodynamic force on the boundary
(represented by $\fG'$) grows faster than the thermodynamic force in the bulk (represented by $f'$). Moreover,
it is easily seen that (A3) is fulfilled for, e.g., the logarithmic case \eqref{flog}. 
\Erem

\gianni{%
\Brem
\label{Generaltrace}
We point out that $\Uad$ actually is a closed, convex, and bounded subset of~$\calX$.
However, it is closed in other spaces as well.
For the reader's convenience, we spend some words on this point.
For $w\in(\Lx2)^3$ with $\div w\in\Lx2$,
the trace $(w\cdot\nu)_{|_\Gamma}$ is a well-defined element of $\HxG{-1/2}$
(in~particular, the definitions \eqref{defcalX} and \eqref{defZ} of $\calX$ and $Z$ are meaningful).
Moreover, the usual integration-by-parts formula holds true in a generalized form
for $w\in(\Lx2)^3$ with $\div w\in\Lx2$ and $v\in\Hx1$.
Namely, we have that $\iO w\cdot\nabla v=-\iO (\div w) v + \<(w\cdot\nu)_{|_\Gamma},v_{|_\Gamma}>$,
where $\<\cpto,\cpto>$ is the duality pairing between $\HxG{-1/2}$ and $\HxG{1/2}$.
In particular, requiring that an element $w\in(\Ldue)^3$ belongs to~$Z$
(i.e.,~it satisfies the conditions $\div w=0$ in $\Omega$ and $w\cdot\nu=0$ on~$\Gamma$)
is~the same as requiring that $\iO w\cdot\nabla v=0$ for every $v\in\Huno$.
Therefore, the whole space $\calX$ can be redefined as the space of $u\in(\LQ\infty\cap\H1\Ldue)^3$
such that $\intQ u\cdot\nabla v=0$ for every $v\in\L2\Huno$.
It follows that $\calX$ is a closed subspace of the Banach space $\tilde\calX:=(\LQ\infty\cap\H1\Ldue)^3$
(as~well as of each of the spaces $(\LQ2)^3$ and $(\H1\Ldue)^3$)
and $\Uad$ is a closed subset of~$\tilde\calX$.
We also notice that, by the above integration-by-parts formula and the assumptions on~$u$, 
we can write the convective term in the next variational formulation
as it is presented in~\eqref{prima} 
(i.e.,~the third integral, to be compared with the second term of~\eqref{ss1}).
\Erem
}%

We now quote some results for the state system \eqref{ss1}--\eqref{ss5} that have recently been proved
in~\cite{CGS13}. 
Prior to this, we notice that the variational form of \eqref{ss1}--\eqref{ss5} reads as follows:  
find functions $(\Mu,\Rho)$ such that
\begin{align}
\label{prima}
  & \iO \dt\rho \, v
  + \iG \dt\rhoG \, \vG
  - \iO \rho u \cdot \nabla v
  + \iO \nabla\mu \cdot \nabla v
  + \iG \nablaG\muG \cdot \nablaG\vG
  = 0
  \non
  \\[1mm]
  & \quad \hbox{\aet\ and for every $(v,\vG)\in\calV$},
   \\[2mm]
   \separa
  \label{seconda}
  & \tauO \iO \dt\rho \, v
  + \tauG \iG \dt\rhoG \, \vG
  + \iO \nabla\rho \cdot \nabla v
  + \iG \nablaG\rhoG \cdot \nablaG\vG
  \non
  \\[1mm]
  & \quad {}
  + \iO f'(\rho) v
  + \iG \fG'(\rhoG) \vG
  = \iO \mu v 
  + \iG \muG \vG
  \non
  \\[1mm]
  & \quad \hbox{\aet\ and for every $(v,\vG)\in\calV$},
     \\[2mm]
  & \rho(0) = \rhoz
  \quad \aeO .
  \label{cauchy}
\end{align}
\Accorpa\State prima cauchy

\noindent
The following result is a combination of the Theorems 2.6, 2.7 and 2.9 in~\cite{CGS13}.

\Bthm
\label{WPstate}
Suppose that the assumptions {\rm (A1)--(A3)} and {\rm (A6)} hold true. 
Then the state 
system {\rm \gian\State} has for every $u\in{\cal U}_R$ a unique solution
$(\Mu,\Rho)$ such that
\begin{equation}
\label{regmurho}
\Mu\in L^\infty(0,T;\calW), \quad \Rho\in W^{1,\infty}(0,T;\calH)\cap H^1(0,T;\calV)\cap L^\infty
(0,T;\calW).
 \end{equation}
Moreover, there are constants $\rho_*,\rho^*\in (-1,1)$ and $\,K_1>0, \,K_2>0$, which depend only on the data
of the state system and $R$, such that the following
holds true:

\noindent {\rm (i)} \quad Whenever $(\Mu,\Rho)$ is the solution to the state system associated with some $u\in{\cal U}_R$,
then we have
\begin{align}
\label{separ}
&\rho_*\le\rho(x,t)\le\rho^* \quad\forall\,(x,t)\in \overline Q\,,\\[2mm]
\label{ssb1}
&\|\Mu\|_{L^\infty(0,T;\calW)}\,+\,\|\Rho\|_{W^{1,\infty}(0,T;\calH)\cap H^1(0,T;\calV)\cap L^\infty
(0,T;\calW)}\,\le\,K_1\,.
\end{align}

\noindent {\rm (ii)} \quad Whenever $u^1,u^2\in {\cal U}_R$ are given and $((\mu^i,\muG^i),
(\rho^i,\rhoG^i))$, $i=1,2$, are the solutions to the corresponding state systems, then 
\begin{align}
\label{stabu}
&\|(\mu^1-\mu^2,\muG^1-\muG^2)\|_{L^\infty(0,T;\calW)}
\non \\
&+\,\|(\rho^1-\rho^2,\rhoG^1-\rhoG^2)\|_{W^{1,\infty}(0,T;\calH)\pier{{}\cap H^1(0,T;\calV){}}\cap L^\infty(0,T;\calW)}\non\\[1mm]
&\le\,K_2\,\|u^1-u^2\|_{H^1(0,T;L^3(\Omega)^3)}\,.
\end{align}
\Ethm

\vspace{3mm}
\Brem
Notice that the pointwise condition \eqref{separ} is meaningful, since it follows from \cite[Sect.~8,~Cor.~4]{Simon}
and \eqref{regmurho} that $\rho\in C^0(\overline Q)$ (and thus, in particular, that $\rhoG\in C^0(\overline\Sigma)$).
\Erem

\pier{We point out that} the uniform separation property 
\eqref{separ} also ensures that the possible singularity 
encoded in the condition  \eqref{singul} never becomes active. 
This implies, in particular, that we may without loss
of generality assume that
\begin{align}
\label{ssb2}
\max_{1\le j\le 3}\left(\|f^{(j)}(\rho)\|_{C^0(\overline Q)}\,+\,\|\fG^{(j)}(\rhoG)\|_{C^0(\overline\Sigma)}
\right)\,\le\,K_1\,,
\end{align}
whenever $\Rho$ is the second component pair of a solution to the state system associated with some $u\in{\cal U}_R$.

\Brem
By virtue of \gian{the well-posedness result given by Theorem~\ref{WPstate}}, 
the control-to-state operator $\calS: u\mapsto (\Mu,\Rho)$ is well defined as a 
mapping between ${\cal U}_R\subset\calX$ and the space defined by the regularity stated in \eqref{regmurho}. 
Moreover, it is Lipschitz
continuous as a mapping from ${\cal U}_R$ into the space 
$$L^\infty(0,T;\calW)\times (W^{1,\infty}(0,T;\calH)\pier{{}\cap H^1(0,T;\calV){}}\cap L^\infty(0,T;\calW)).$$
\Erem


\section{Fr\'echet differentiability of the control-to-state operator}
\label{FRECHET}
\setcounter{equation}{0}
In this section, we aim to show the Fr\'echet differentiability of the control-to-state 
operator $\calS$ in suitable Banach spaces. Throughout this section, we assume that
$\ub\in {\cal U}_R$ is fixed and that the general assumptions (A1)--(A3), (A5) and (A6) are satisfied,
so that the global estimates \eqref{ssb1} and \eqref{ssb2} are valid for the associated solution 
$(\mub,\muGb,\rhob,\rhoGb)=\calS(\ub)$ to the state system. We then consider the
linearized system, \gian{where $h\in\calX$,} 
\Bsist
  && \iO \dt\xi \, v
  + \iG \dt\xiG \, \vG
  + \iO \nabla\eta \cdot \nabla v
  + \iG \nablaG\etaG \cdot \nablaG\vG
  \non
  \\
  && \quad {}
  - \iO \xi \ub \cdot \nabla v
  - \iO \rhob \, h \cdot \nabla v
  = 0
  \non
  \\
  && \quad \hbox{\aet\ and for every $(v,\vG)\in\calV$},
  \label{primaL}
  \\
  \separa
  && \tauO \iO \dt\xi \, v
  + \tauG \iG \dt\xiG \, \vG
  + \iO \nabla\xi \cdot \nabla v
  + \iG \nablaG\xiG \cdot \nablaG\vG
  \non
  \\
  && \quad {}
  + \iO f''(\rhob) \xi v
  + \iG f''_\Gamma(\rhoGb)\xiG \vG
  = \iO \eta v 
  + \iG \etaG \vG
  \non
  \\
  && \quad \hbox{\aet\ and for every $(v,\vG)\in\calV$},
    \label{secondaL}
  \\
  && \xi(0) = 0  \quad \aeO, \quad\xiG(0)=0 \quad\aeG.
  \label{cauchyL}
\Esist
\Accorpa\Linear primaL cauchyL
\gianni{Since $\ub\in{\calU}_R\subset\calX$ and $h\in\calX$, then $\div\ub=0$ and $\ub\cdot\nu=0$,
and the same conditions hold for~$h$,
so that \Linear\ is the weak form of the linear initial-boundary value problem}
\begin{align}
\label{L1}
& \gianni{\dt\xi - \Delta\eta = \pier{{}- \nabla\xi\cdot\ub -\nabla\rhob\cdot h {}}}
\quad\aeQ,
\\[1mm]
\label{L2}
& \gianni{\dt\xiG + \dn\eta - \DeltaG\etaG = 0 }
\aand
\eta_{|\Sigma}=\etaG, \,\quad\aeS,\\[1mm]
\label{L3} 
&\tauO\dt\xi-\Delta\xi+f''(\rhob)\xi=\eta \quad\aeQ,\\[1mm]
\label{L4}
&\tauG\dt\xiG+\dn\xi-\DeltaG\xiG+f''_\Gamma(\rhoGb)\xiG=\etaG 
\quad\,\mbox{and}\,\quad\xi_{|\Sigma}=\xiG,\,\quad\aeS,\\[1mm]
\label{L5}
&\xi(0)=0 \quad\aeO,\qquad \xiG(0)=0 \quad
\aeG.
\end{align}
\gianni{However, we only refer to the problem in the form \Linear.}

We expect the following to hold true: 
if the system \Linear\ admits for every \gian{$h\in\calX$} a unique solution 
$((\eta,\etaG),(\xi,\xiG))=:((\eta^h,\etaG^h),(\xi^h,\xiG^h))$ in a
suitable Banach space, then the Fr\'echet derivative $D\calS(\ub)$ of $\calS$ 
at $\ub$ (if it exists), evaluated at $h$, should have the form 
$D\calS(\ub)(h)=((\eta^h,\etaG^h),(\xi^h,\xiG^h))$.

\Bthm
Suppose that the assumptions {\rm (A1)--(A3), (A5)} and {\rm (A6)} are fulfilled, let $\ub\in{\cal U}_R$
be given, and let $((\mub,\muGb),(\rhob,\rhoGb))=\calS(\ub)$ be the
associated unique solution to the state system \eqref{ss1}--\eqref{ss5} having
the regularity properties stated in \eqref{regmurho}. Then the system \Linear\ has for every
$h\in{\cal X}$ a unique solution $((\eta,\etaG),(\xi,\xiG))$ such that
\begin{align}
\label{regeta}
(\eta,\etaG)\in L^2(0,T;\calW),\ \quad (\xi,\xiG)\in H^1(0,T;\calH)\cap L^\infty(0,T;\calV)
\cap L^2(0,T;\calW). 
\end{align}
Moreover, the linear mapping $h\mapsto ((\eta,\etaG),(\xi,\xiG))$ is continuous as a mapping
from ${\cal X}$ into the space
\begin{equation}\label{defcalY}
{\cal Y}:=L^2(0,T;\calV)\times (H^1(0,T;\calH)\cap L^\infty(0,T;\calV))\,.
\end{equation}
\Ethm

\proof
We employ a slightly modified 
Faedo-Galerkin scheme with a proper choice of the Hilbert basis. To this end,
we introduce the operator $\calA\in\calL(\calV;\calVp)$ by setting
\Beq
  \< \calA(w,\wG) , (v,\vG) >_{\calV}
  := \iO \nabla w \cdot \nabla v
  + \iG \nablaG\wG \cdot \nablaG\vG
  \quad \hbox{for $(w,\wG),(v,\vG)\in\calV$},
  \label{defA}
\Eeq
and notice that $\calA$ is nonnegative and weakly coercive.
Indeed, we have that
\Beq
  \< \calA(v,\vG) , (v,\vG) >_{\calV}
  + \normaHH{(v,\vG)}^2 = \normaVV{(v,\vG)}^2 
  \quad \hbox{for every $(v,\vG)\in\calV$}.
  \label{coerc}
\Eeq
Moreover, as the embedding $\calV\subset\calH$ is compact,
the resolvent of $\calA$ is compact as well, and
the spectrum of $\calA$ reduces to a discrete set of eigenvalues,
the eigenvalue problem being
\Beq
  (e,e_\Gamma) \in \calV \setminus \{(0,0)\}
  \aand
  \calA(e,e_\Gamma) = \lambda(e,e_\Gamma) \,.
  \label{eigenpbl}
\Eeq
More precisely, we can rearrange the eigenvalues and choose the eigenvectors in such a way that
\begin{align}
  & 0 = \lambda_1 < \lambda_2 \leq \lambda_3 \leq \dots
  \aand
  \lim_{j\to\infty} \lambda_j = +\infty ,
  \label{eigenvalues}
  \\
  & \calA \Ej = \lambda_j \Ej
  \aand\non\\
  &\quad \iO \ei\ej + \iG \eGi\eGj = \delta_{ij} \pier{{}:= \begin{cases}
  \, 1 &\hbox{if } i=j\\
  \, 0 &\hbox{if } i\not= j
  \end{cases}}  \quad \hbox{for $i,j=1,2,\dots$},
  \qquad
  \label{eigenvectors}
\end{align}
and $\graffe{\Ej}$ generates a dense subspace of both $\calV$ and  $\calH$.
We notice that
\Beq
  \iO \nabla\ei \cdot \nabla\ej
  + \iG \nablaG\eGi \cdot \nablaG\eGj
  = \lambda_i \Bigl( \iO \ei\ej + \iG \eGi\eGj \Bigr)
  = \lambda_i \delta_{ij} 
  \quad \hbox{for $i,j=1,2,\dots$}.
  \non
\Eeq
We also observe that every element $(w,\wG)\in\calH$ can be written~as
\Beq
  (w,\wG) = \somma j1\infty w_j \Ej
  \quad \hbox{with} \quad
  \somma j1\infty |w_j|^2 = \normaHH{(w,\wG)}^2 < + \infty,
  \non
\Eeq
and that (on account of \eqref{coerc})
\Beq
  (w,\wG) \in \calV
  \quad \hbox{if and only if} \quad
  \somma j1\infty (1+\lambda_j) |w_j|^2 < + \infty \,.
  \non
\Eeq
Namely, the last sum yields the square of a norm on $\calV$
that is equivalent to~$\normaVV\cpto$. 

At this point, we set
\Beq
  \calVn := \Span\graffe{\Ej:\ 1\leq j\leq n}
  \aand
  \calV_\infty := \bigcup_{j=1}^\infty \calVn
  = \Span\graffe{\Ej:\ j\geq1},
  \label{defVn}
\Eeq
and, for every~$n\geq 1$, we look for a quadruple $(\etan,\etaGn,\xin,\xiGn)$ satisfying
\Bsist
  && \Etan \in \L2\calVn
  \aand
  \Xin \in \H1\calVn,
  \vphantom\sum
  \label{regsoluzn}
  \\[1mm]
\label{priman}
  && \iO \dt\xin \, v + \iG \dt\xiGn \, \vG
  - \iO \xin \ub \cdot \nabla v-\iO \rhob h \cdot\nabla v 
  + \iO \nabla\etan \cdot \nabla v
  + \iG \nabla\etaGn \cdot \nabla\vG\non\\[1mm]
  && \quad + \frac 1n \iO \etan v + \frac 1n \iG \etaGn \vG 
  = 0 \quad \hbox{\aet\ and for every $(v,\vG)\in\calVn$},
 \\[2mm]
 \label{secondan}
  && \tauO \iO \dt\xin \, v
  + \tauG \iG \dt\xiGn \, \vG
  + \iO \nabla\xin \cdot \nabla v
  + \iG \nablaG\xiGn \cdot \nablaG\vG
  \non
  \\[1mm]
  && \quad {}
  + \iO f''(\rhob)\xin v + \iG \fG''(\rhoGb)\xiGn \vG
  = \iO \etan v
  + \iG \etaGn \vG
  \non
  \\[1mm]
  &&\qquad \hbox{\aet\ and for every $(v,\vG)\in\calVn$},
    \\[2mm]
    \label{inin}
  && \xin(0) = 0
  \quad \aeO, \qquad  \pier{\xin_\Gamma(0) = 0
  \quad \aeG.}  \label{cauchyn}
\Esist
\Accorpa\Pbln regsoluzn cauchyn

\step
Existence for the discrete problem

For every fixed $n\in\enne$, we are looking for $\Etan$ and $\Xin$ in the form
\Beq
  \Etan(t) = \somma j1n \eta_j^n(t) \Ej
  \aand
  \Xin(t) = \somma j1n \xi_j^n(t) \Ej\,,
  \non
\Eeq
for some $\eta_j^n\in L^2(0,T)$ and $\xi_j^n\in H^1(0,T)$, $1\le j\le n$. 
Let us introduce the $n$-vector functions $\overline\eta:=(\eta_1^n,\ldots, \eta_n^n)$ and 
$\overline\xi:=(\xi_1^n,\ldots,\xi_n^n)$.
Then, making the special choices $(v,\vG)=\Ei$ for $i=1,\dots,n$, we can rewrite the
system \accorpa{priman}{secondan} in the form
\Bsist
  && \pier{{\overline\xi\,}'(t) -U(t) \, \overline\xi(t) + D_n \, 
  \overline\eta(t)} = C(t)
  \non
  \\
  && \aand
  B \, {\overline\xi\,}'(t) + D \, \overline\xi(t) + \pier{{}G(t) \, \overline\xi(t)} = \overline\eta(t),
  \label{odesystem}
\Esist
where 
$D_n:=\diag(\lambda_1+\frac 1n,\dots,\lambda_n+\frac 1n)$,
$D:=\diag(\lambda_1,\dots,\lambda_n)$ 
and where the matrices $U=(u_{ij})\in\L2{\erre^{n\times n}}$, 
$G=(g_{ij})\in\L\infty{\erre^{n\times n}}$, and $B=(b_{ij})\in\erre^{n\times n}$ 
and the vector $C=(c_i)\in L^2(0,T;\erren)$ are given~by
\begin{align}
 & u_{ij}(t) := \iO \ej \ub(t) \cdot \nabla\ei, \quad 
 \pier{g_{ij}(t) := \iO f''(\rhob)(t) \ej \ei + \iG \fG''(\rhoGb)(t) \eGj \eGi,} \non\\[1mm]
 &   b_{ij} := \tauO \iO \ej \ei + \tauG \iG \eGj \eGi,
 \quad
 \pier{c_i(t):=\iO \rhob(t)\,h(t)\cdot\nabla e^i,}  
 \quad\mbox{for $i,j=1,\dots,n$}.
  \non
\end{align}

By adding the second identity in 
\eqref{odesystem} to the first one multiplied by~$D_n^{-1}$,
we obtain the equivalent system
\begin{align*}
  (D_n^{-1}+B) \, {\overline \xi\,}'(t) + V(t) \, \pier{\overline\xi(t)} = D_n^{-1}C(t)
  \aand
  \overline\eta(t) = B \, {\overline\xi\,}'(t) + D \, \overline\xi(t) + \pier{{}G(t) \, \overline\xi(t)},
\end{align*}
where $V:=D \pier{{}+ G} -D_n^{-1}U$ belongs to $\L2{\erre^{n\times n}}$ 
and $D_n^{-1}+B$ is invertible, as we now verify.
To this end, we show that $B$ is positive definite.
Indeed, for any vector $y=(y_1,\dots,y_n)\in\erre^n$, by setting $(v,\vG):=\somma j1n y_j\Ej$,
we have that
\Bsist
  && (By) \cdot y = \somma {i,j}1n b_{ij} y_j y_i
  \,=\, \tauO \iO \somma i1n y_i \ei \somma j1n y_j \ej
  + \tauG \iO \somma i1n y_i \eGi \somma j1n y_j \eGj
  \non
  \\
  && =\, \tauO \iO |v|^2 + \tauG \iG |\vG|^2
  \,\geq\, \min\{\tauO,\tauG\}\, \normaHH{(v,\vG)}^2
  \,=\, \min\{\tauO,\tauG\}\, \norma y_{\erre^n}^2 \,.
  \non
\Esist
Hence, $D_n^{-1}+B$ is positive definite as well,  and thus invertible.
Therefore, by virtue \gian{of} standard results for initial value problems for ordinary differential
equations, the discrete problem \Pbln\ has a unique solution having the 
asserted regularity.

At this point, our aim is to show that the solutions to the \juerg{discrete} problem
converge to a solution to \eqref{primaL}--\eqref{cauchyL} as $n$ tends to infinity, 
at least for a subsequence.
To this end, we start estimating and find bounds that do not depend on~$n$. In the following,
$C_i$, $i\in\enne$, will denote positive constants that may depend on the data of the system
and on~$R$, but not on $n\in\enne$.

\step
First a priori estimate

We test \eqref{priman}, written at the time $s$, by $\Etan(s)$ and integrate over $(0,t)$ with respect to~$s$ to find that
\Bsist
  && \intQt \dt\xin \, \etan
  + \intSt \dt\xiGn \, \etaGn
  + \intQt |\nabla\etan|^2
  + \intSt |\nablaG\etaGn|^2
  \non
  \\
  && \quad {}
  + \frac 1n \intQt |\etan|^2
  + \frac 1n \intSt |\etaGn|^2
  = \intQt \xin \ub \cdot \nabla\etan + \intQt \rhob\,h\cdot \nabla\etan \,.
  \non
\Esist
Next, we test \eqref{secondan} by $\dt\Xin(s)$, integrate over $(0,t)$ with respect to~$s$,
and add the expression \,$\intQt\xin\dt\xin + \intSt\xiGn\dt\xiGn$ to both sides, for convenience.
We \pier{infer} that
\begin{align}
  & \tauO \intQt |\dt\xin|^2
  + \tauG \intSt |\dt\xiGn|^2
  + \frac 12 \, \normaVV{\Xin(t)}^2\non\\[1mm]
  &= \intQt \left(1-f''(\rhob)\right)\xin \dt\xin
  + \intSt \left(1-\fG''(\rhob)\right)\xiGn \dt \xiGn
    + \intQt \etan \dt\xin
  + \intSt \etaGn \dt\xiGn\,. 
    \non
\end{align}
At this point, we add these equalities
and notice that four terms cancel and that the remaining terms on the \lhs\ are nonnegative.
Moreover, we use the global estimates \eqref{ssb1}, \eqref{ssb2}, and Young's inequality.
We obtain  that
\Bsist
  && \intQt |\nabla\etan|^2
  + \intSt |\nablaG\etaGn|^2
  + \tauO \intQt |\dt\xin|^2
  + \tauG \intSt |\dt\xiGn|^2
  + \frac 12 \, \normaVV{\Xin(t)}^2
  \non
  \\
  && \leq \intQt |\xin| \, |\ub| \, |\nabla\etan|
  + \intQt |\rhob|\,|h|\,|\nabla\etan|
  + \frac {\tauO}2 \intQt |\dt\xin|^2
  + \frac {\tauG} 2 \intSt |\dt\xiGn|^2\non\\
  && \quad {}
  + C_1 \intQt |\xin|^2
  + C_1 \intSt |\xiGn|^2 \,.
  \non
\Esist
On the other hand, the \Holder, Sobolev and Young inequalities yield that
\Bsist
  && \intQt |\xin| \, |\ub| \, |\nabla\etan|
  \leq \iot\|\xin(s)\|_6 \, \norma{\ub(s)}_3 \, \norma{\nabla\etan(s)}_2 \, ds
  \non
  \\
  && \leq \frac 14 \intQt|\nabla\etan|^2
  + C_2 \, \pier{{}\norma{\ub}_{\L\infty {L^3(\Omega)}}^2} \iot  \norma{\xin(s)}_V^2 \, ds\,.
  \non
\Esist
\pier{Moreover, we have that}
\begin{align*}
&\intQt |\rhob|\,|h|\,|\nabla\etan|\,\le\,\iot\|\rhob(s)\|_\infty\|h(s)\|_2\,\|\nabla\etan(s)\|_2\,ds
\non\\
&\le\, \frac 14\intQt|\nabla\etan|^2+C_3\intQt|h|^2\,.
\end{align*}
Therefore, rearranging and applying Gronwall's lemma, we can infer that for all $t\in (0,T]$ it holds 
\Beq
\label{einsn}
  \pier{\bigg(\intQt|\nabla\etan|^2 + \intSt|\nablaG \etaGn|^2 \bigg)^{\!\!1/2}}
    + \norma\Xin_{H^1(0,t;\calH)\cap L^\infty(0,t;\calV)}
  \leq C_4\|h\|_{L^2(0,t;H)^3}\,.
\Eeq

\step
Second a priori estimate

We insert $(v,\vG)=(\etan,\etaGn)$ in \eqref{secondan}. As it turns out, all of the
resulting terms can be handled directly by means of Young's inequality to yield an inequality
of the form
\begin{align*}
&\intQt |\etan|^2+\intSt |\etaGn|^2 \,\leq \,C_5\Big(\|(\xin,\xiGn)\|_{H^1(0,T;\calH)
\cap L^2(0,T;\calV)}^2+\intQt|\nabla\etan|^2+\intSt|\nablaG \etaGn|^2\Big),
\end{align*}
and it follows from \eqref{einsn} that 
\begin{equation}
\label{zwein}
\|\Etan\|_{L^2(0,t;\calV)}\,\le\,C_6\|h\|_{L^2(0,t;H)^3}\,.
\end{equation}

\step
Existence of a unique solution to the linearized system

We account for \accorpa{einsn}{zwein} 
and use standard weak and weak star compactness results,
as well as \cite[Sect.~8, Cor.~4]{Simon}).
It follows that, as $n$ tends to infinity,
\Bsist
  & \Etan \to (\eta,\etaG)
  & \quad \hbox{weakly in $\L2\calV$},
  \label{convetan}
  \\
  & \pier{\frac1n \Etan \to (0,0)}
  & \quad \pier{\hbox{strongly in $\L2\calV$},}
  \label{convetanbis}
  \\
& \Xin \to \Xi
  & \quad \hbox{weakly star in $\H1\calH\cap L^\infty(0,T;\calV)$}
  \non
  \\
  && \quad \hbox{and strongly in $C^0([0,T];\calH)$},
  \label{convxin}
\Esist
at least for a subsequence, which is again indexed by $n$.
In particular, we have $\xi(0)=0$ and $\xiG(0)=0$.
We also recall that $\ub,f''(\rhob) \in L^\infty(Q)$ and 
$\fG''(\rhoGb)\in L^\infty(\Sigma)$, which implies~that
\begin{align*}
&\xin\ub\to \xi\ub \quad\mbox{strongly in }\,L^2(0,T;H)^3\,,\\
&f''(\rhob)\xin\to f''(\rhob)\xi\quad\mbox{strongly in }\,L^2(Q)\,,\\
&\fG''(\rhoGb)\xiGn\to \fG''(\rhoGb)\xiG \quad\mbox{strongly in } L^2(\Sigma)\,.
\end{align*}
Now, we recall \eqref{defVn} for the definition of $\calV_\infty$,
and take an arbitrary $\calV_\infty$-valued step function~$(v,\vG)$.
Since the range of $(v,\vG)$ is finite-dimensional,
there exists some\, $m\in\enne$\, such that
$(v,\vG)(t)\in\calV_m$ \aat.
It follows that $(v,\vG)(t)\in\calVn$ \aat\ and every $n\geq m$,
so that we can test \eqref{priman} and \eqref{secondan}, written at the time~$t$,
 by $(v,\vG)(t)$
and integrate over~$(0,T)$.
At this point, it is \sfw\ to deduce that
$(\eta,\etaG)$ and $\Xi$ satisfy the integrated 
version of \eqref{primaL}--\eqref{cauchyL} for every such step \juerg{function}, namely, we have that
\Bsist
  && \intQ \dt\xi \, v + \intS \dt\xiG \, \vG
  - \intQ \xi \ub \cdot \nabla v -\intQ \rhob h \cdot\nabla v
  + \intQ \nabla\eta \cdot \nabla v
  + \intS \nabla\etaG \cdot \nabla\vG
  = 0\,,
  \quad
  \non
  \\
  \separa
  && \tauO \intQ \dt\xi \, v
  + \tauG \intS \dt\xi \, \vG
  + \intQ \nabla\xi \cdot \nabla v
  + \intS \nablaG\xiG \cdot \nablaG\vG
  \non
  \\
  && \quad {}
  + \intQ f''(\rhob)\xi  v
  + \intS \fG''(\rhoGb)\xiG  \vG
  = \intQ \eta v 
  + \intS \etaG \vG \,.
  \non
\Esist
By density, the same equations hold true for every $(v,\vG)\in\L2\calV$.
This implies that \eqref{primaL}--\eqref{secondaL} hold true \aet\ and for every $(v,\vG)\in\calV$,
as desired. It is thus shown that $((\eta,\etaG),\Xi)$ is a solution to the
linearized system \Linear.

Next, we show that there can be no  other such solution. To this end, assume that $((\eta^i,\etaG^i),
(\xi^i,\xiG^i))$,
$i=1,2$, are two solutions such that 
$$
(\eta^i, \etaG^i)\in L^2(0,T;\calV) \quad\mbox{and}\quad (\xi^i,\xiG^i)\in H^1(0,T;\calH)\cap
L^\infty(0,T;\calV) \quad\mbox{for }\, i=1,2.
$$
We put $(\eta,\etaG):=(\eta^1,\etaG^1)-(\eta^2,\etaG^2)$ and $(\xi,\xiG):=(\xi^1,\xiG^1)-(\xi^2,\xiG^2)$.
Then $((\eta,\etaG),(\xi,\xiG))$ solves the system \Linear, where, in this case, the expression 
\,$\iO \rhob h\cdot \nabla v$\, is not present. Now, we repeat the two a priori estimates performed
above for the approximating system, but this time we \pier{proceed directly on} the system \Linear. We then
recover the estimates \eqref{einsn} and \eqref{zwein}, but this time with zero right-hand sides.
Hence, $(\eta,\etaG)$ and $(\xi,\xiG)$ vanish, which proves the uniqueness.
   
\step
Further regularity

We still need to show that $(\xi,\xiG)\in L^2(0,T;\calW)$. This is an immediate consequence of
\cite[Lem.~3.1]{CGS13}: indeed, we can write \eqref{secondaL} in the form
\begin{equation}
\iO \nabla\xi\cdot\nabla v+\iG \nablaG\xiG\cdot\nablaG\vG=\iO g v+\iG g_\Gamma \vG 
\quad \forall (v,\vG)\in\calV,
\label{pier1}
\end{equation} 
for a.e. $t\in (0,T)$, where we define $(g,g_\Gamma)\in L^2(0,T;\calH)$ by
$$
g:=\eta-f''(\rhob)\xi-\tauO \dt\xi\,,\quad g_\Gamma:=\etaG-\fG''(\rhoGb)\xiG-\tauG \dt\xiG\,.
$$
Obviously, $(g(t),g_\Gamma(t))\in\calH$ for a.e. $t\in(0,T)$. It then follows from \cite[Lem.~3.1]{CGS13}
that, for a.e. $t\in (0,T)$, it holds $(\xi(t),\xiG(t))\in\calW$, as well as 
\Beq
\|(\xi(t),\xiG(t))\|_{\calW}\,\le\,C_\Omega\left(\|(\xi(t),\xiG(t))\|_{\calV}\,+\,\|(g(t),g_\Gamma(t))\|_\calH
\right), \label{pier2} 
\Eeq
with a constant $C_\Omega>0$ that depends only on $\Omega$. Since we have $(\xi,\xiG)\in L^\infty(0,T;\calV)$,
we conclude that indeed $(\xi,\xiG)\in L^2(0,T;\calW)$. \pier{Arguing as above on 
the equation \eqref{primaL}, to be written similarly as in \eqref{pier1}, and observing that 
\begin{align*}
&\norma{\nabla\xi\cdot\ub -\nabla\rhob\cdot h}_{\L2H} \\
&\leq C_7 \Big( \norma{\nabla\xi}_{\L2{L^6(\Omega)^3}} \norma{\ub}_{\L\infty{L^3(\Omega)^3}} + 
\norma{\nabla\rhob}_{\L\infty{L^6(\Omega)^3}} \norma{h}_{\L2{L^3(\Omega)^3}} \Big) \\
&\leq C_8 \Big( \norma{(\xi,\xiG)}_{\L2{\calW}} \norma{\ub}_{\L\infty{L^3(\Omega)^3}} + 
\|(\rhob, \rhob_\Gamma) \|_{L^\infty (0,T;\calW)} \norma{h}_{\L2{L^3(\Omega)^3}} \Big)\,, 
\end{align*} 
it is not difficult to conclude that $(\eta,\etaG)\in L^2(0,T;\calW)$ (cf.~\cite[Sect.~3]{CGS13}).}

\pier{At this point, it} remains to show the asserted continuity properties of the mapping $\,h\mapsto ((\eta,\etaG),\Xi)$.
Now, it follows from the weak and weak star sequential semicontinuity of norms and from the estimates \eqref{einsn} and  \eqref{zwein} that, for every $h\in\calX$,  
\begin{equation}
\label{lincon}
\|(\eta,\etaG)\|_{L^2(0,T;\calV)}\,+\,\|\Xi\|_{H^1(0,T;\calH)
\cap L^\infty(0,T;\calV)} \,\le \, \pier{C_9\,\|h\|_{\cal X}}\,.
\end{equation}
The  assertion is thus completely proved.
\QED

\vspace{3mm}
We now turn our interest to the Fr\'echet differentiability. 
\gianni{We recall the definitions \eqref{defcalX} and \eqref{defcalY} of the spaces $\calX$ and $\calY$
and prove the following result.}
\Bthm
\label{DcalS}
Assume that {\rm (A1)--(A3), (A5)} and {\rm (A6)} are fulfilled. 
Then the control-to-state operator $\calS$ is Fr\'echet
differentiable at every $\ub\in {\cal U}_R$ as a mapping from the space
\gianni{${\cal X}$} into the space~${\cal Y}$. 
Moreover, for every 
$\ub\in {\cal U}_R$ and every $h\in {\cal X}$ we have that the
Fr\'echet derivative $D\calS(\ub)$ of $\calS$ at $\ub$ satisfies
$D\calS(\ub)(h)= (\eta,\etaG,\xi,\xiG)$, which is the unique solution
to the linearized system \Linear\ associated with $h$. 
\Ethm

\proof 
Since ${\cal U}_R$ is open, there is some $\Lambda>0$ such that
$\ub+h\in {\cal U}_R$ whenever $h\in{\cal X}$ and $\|h\|_{\cal X}\le\Lambda$. In the
following, we consider only such perturbations $h$, for which we define the quantities
\begin{align*}
&((\mu^h,\muG^h),(\rho^h, \rhoG^h)):=\calS(\ub+h), \quad y^h:=\rho^h-\rhob-\xi^h,
\quad y^h_\Gamma:=\rhoG^h-\rhoGb-\xiG^h,\\[1mm]
&z^h:=\mu^h-\mub-\eta^h,\quad z^h_\Gamma:=\muG^h-\muGb-\etaG^h\,.
\end{align*} 
Obviously, we have $y^h_{|\Sigma}=y^h_\Gamma$ and $z^h_{|\Sigma}=z^h_\Gamma$, as well~as
\begin{align}
\label{regyzh}
&(y^h,y^h_\Gamma)\in H^1(0,T;\calH)\cap L^\infty(0,T;\calV)\,,\quad
(z^h,z^h_\Gamma)\in L^2(0,T;\calV)\,.   
\end{align} 
Since we know already from the previous theorem that the linear mapping 
$h\mapsto ((\eta^h,\eta^h_\Gamma),\linebreak
(\xi^h,\xiG^h))$ is continuous as a mapping from ${\cal X}$ into ${\cal Y}$, it
suffices to show that there is an increasing mapping $Z:(0,\Lambda)\to (0,+\infty)$ 
such that $\,\lim_{\lambda\searrow 0} Z(\lambda)/\lambda^2=0\,$ and
\begin{equation}
\label{Frechet}
\|(z^h,z^h_\Gamma)\|^2_{L^2(0,T;\calV)}+\|(y^h,y^h_\Gamma)\|_{H^1(0,T;\calH)\cap
L^\infty(0,T;\calV)}^2\,\le\,Z(\|h\|_{\cal X})\,.
\end{equation}
We also recall that $((\mu^h,\muG^h),(\rho^h, \rhoG^h))$ satisfy the global estimates
stated in \eqref{ssb1}, \eqref{ssb2}, and we observe that it follows from  Taylor's 
theorem that there is some $C_1>0$ such that 
\begin{align}
\label{taylor1}
& \left|f'(\rho^h)-f'(\rhob)-f''(\rhob)\xi^h\right|\,\le\,C_1\left|y^h\right|+C_1
\left|\rho^h-\rhob\right|^2 \quad
\aeQ,\\[1mm]
\label{taylor2}
&\left|\fG'(\rhoG^h)-\fG'(\rhoGb)-\fG''(\rhoGb)\xiG^h\right|\,\le\,C_1
\left|y^h_\Gamma\right|+C_1\left|\rhoG^h-\rhoGb\right|^2
\quad\aeS,
\end{align}
where, here and in the remainder of this proof, $C$ and $C_i$, $i\in\enne$, denote positive
constants that may depend on the data of the system and $R$, but not on the special
choice of $h$ with $\|h\|_\calX\le\Lambda$. Moreover, using the state equations and
the linearized system, we readily verify that the following identities are valid:
\begin{align}
\label{yhzh1}
&\iO \dt y^h v+\iG \dt y_\Gamma^h \vG+\iO \nabla z^h\cdot\nabla v+\iG \nablaG z^h_\Gamma
\cdot\nablaG \vG=\iO y^h\ub\cdot\nabla v+\iO (\rho^h-\rhob)h\cdot \nabla v\nonumber\\[1mm]
&\quad \mbox{for all }\,(v,\vG)\in \calV \,\quad\mbox{and }\,\aet,\\[2mm]
\label{yhzh2}
&\tauO\iO\dt y^h v+\tauG\iG\dt y^h_\Gamma \vG+\iO\nabla y^h\cdot \nabla v
+\iG\nablaG y^h_\Gamma\cdot\nablaG \vG\nonumber\\[1mm]
&+\iO\left(f'(\rho^h)-f'(\rhob)-f''(\rhob)\xi^h\right)v +\iO\left(\fG'(\rhoG^h)
-\fG'(\rhoGb)-\fG''(\rhoGb)\xiG^h\right)\vG\nonumber\\[1mm]
&=\iO z^h v +\iG z^h_\Gamma \vG \qquad \mbox{for all }\,(v,\vG)\in \calV \,\quad\mbox{and }\,\aet.
\end{align}

\step
First estimate

For $s\in (0,T)$, we insert $(v,\vG)=(z^h(s),z^h_\Gamma(s))$ in \eqref{yhzh1} and
$(v,\vG)=(\dt y^h(s),\dt y^h_\Gamma(s))$ in \eqref{yhzh2}. \pier{The last position is formal, but the following computations can be justified rigorously by arguing, e.g., as in \cite[Appendix]{CGG3}}.  
We then add the two resulting identities\pier{, integrate
over $(0,t)$, where $t\in (0,T)$, and add on both sides the quantity \,$\intQt y^h\,\dt y^h+\intSt y^h_\Gamma\,\dt y^h_\Gamma$. Observing that some terms cancel out and using}
the inequalities \eqref{taylor1} and \eqref{taylor2}, we arrive at the inequality
\begin{align}
\label{diff1}
&\intQt|\nabla z^h|^2+\intSt|\nablaG z^h_\Gamma|^2+\frac 12 \|(y^h(t),y^h(t))\|^2_{\calV}
+\tauO\intQt|\dt y^h|^2+\tauG\intSt|\dt y^h_\Gamma|^2\non\\[1mm]
&\le\,\intQt |y^h||\ub||\nabla z^h| +\intQt |\rho^h-\rhob||h| |\nabla z^h|+C\intQt|y^h||\dt y^h|
+C\intSt|y^h_\Gamma||\dt y^h_\Gamma|\non\\[1mm]
&\quad + C\intQt|\dt y^h||\rho^h-\rhob|^2+C\intSt|\dt y^h_\Gamma||\rhoG^h-\rhoGb|^2
=:\sum_{j=1}^6 I_j,
\end{align}
with obvious notation. We estimate the six terms on the \rhs\ individually, using 
the H\"older, Young and Sobolev inequalities, and invoking \eqref{stabu}. We obtain the
following estimates:
\begin{align}
\label{diff2}
&I_1\,\le\,\int_0^t\|y^h(s)\|_6\,\|\ub(s)\|_3\,\|\nabla z^h(s)\|_2\,ds \non \\
&\quad \le\,\frac 14\intQt |\nabla z^h|^2+C\|\ub\|_{L^\infty(0,t;L^3(\Omega)^3)}^2
\|y^h\|_{L^2(0,t;V)}^2\non\\
&\quad \le\,\frac 14\intQt |\nabla z^h|^2+C\,\|y^h\|_{L^2(0,t;V)}^2\,,
\\[2mm]
\separa
\label{diff3}
&I_2\,\le\,\int_0^t\|\rho^h(s)-\rhob(s)\|_6\,
\|h(s)\|_3\,\|\nabla z^h(s)\|_2\,ds\non\\
&\quad \le\,\frac 14\intQt|\nabla z^h|^2+C\,\|h\|^2_{L^\infty(0,t;L^3(\Omega)^3)}\|\rho^h-\rhob\|_{L^2(0,t;V)}^2\non\\
&\quad \le\,\frac 14\intQt|\nabla z^h|^2+C\,\|h\|_{H^1(0,t;L^3(\Omega)^3)}^4\,,
\\[2mm]
\separa
\label{diff4}
&I_3\,\le\,\frac{\tauO}4\intQt|\dt y^h|^2+C\intQt|y^h|^2\,,
\\[2mm]
\separa
\label{diff5}
&I_4\,\le\,\frac {\tauG}4\intSt|\dt y^h_\Gamma|^2+C\intSt|y^h_\Gamma|^2\,,
\\[2mm]
\label{diff6}
&I_5\,\le\,C\int_0^t\|\dt y^h(s)\|_2\,\|\rho^h(s)-\rhob(s)\|_4^2\,ds \non \\
&\quad \le\,\frac{\tauO}4\intQt|\dt y^h|^2+C\int_0^t\|\pier{(\rho^h-\rhob)(s)}\|_V^4\,ds
\non
\\
&\quad\le\,\frac{\tauO}4\intQt |\dt y^h|^2
+ C\,\|h\|^4_{H^1(0,t;L^3(\Omega)^3)}\,,
\end{align}
and, by the same token,
\begin{align}
\label{diff7}
&I_6\,\le\,\frac{\tauG}4\intSt|\dt y^h_\Gamma|^2+C\,\|h\|^4_{H^1(0,t;L^3
(\Omega)^3)}\,.
\end{align}
At this point, we can combine the estimates 
\eqref{diff1}--\eqref{diff7} and infer from Gronwall's lemma that, for every $t\in (0,T]$,
\gian{\Beq
 \|(y^h,y^h_\Gamma)\|^2_{H^1(0,t;\calH)\cap L^\infty(0,t;\calV)}
\,+\,\norma{(\nabla z^h,\nablaG z^h_\Gamma)}_{L^2(0,t;{\calH^3})}^2
\leq C_2\,\|h\|^4_{H^1(0,t;L^3(\Omega)^3)}.
\label{diff8}
\Eeq
}%

\step Second estimate

\vspace{1mm} \noindent
Next, we insert, for $s\in (0,T)$, $(v,\vG)=(z^h(s),z^h_\Gamma)(s)$ in \eqref{yhzh2}
and integrate the resulting equation over $(0,t)$, where $t\in (0,T]$. Using \eqref{taylor1} and 
\eqref{taylor2} once more, we then arrive at the estimate
\begin{align}
\label{diff9}
&\intQt|z^h|^2+\intSt|z^h_\Gamma|^2\,\le\,\tauO\intQt|\dt y^h||z^h|+
\tauG\intSt|\dt y^h_\Gamma||z^h_\Gamma|+\intQt|\nabla y^h||\nabla z^h|
\non\\
&\quad +\intSt \pier{|\nablaG y^h_\Gamma|\, |\nablaG z^h_\Gamma|}+C\intQt|y^h|\, |z^h|+C\intSt|y^h_\Gamma|
|z^h_\Gamma|+C\intQt|\rho^h-\rhob|^2|z^h|\non\\
&\quad + C\intSt|\rhoG^h-\rhoGb|^2|z^h_\Gamma|\,.
\end{align}
The sum of the first six summands on the \rhs, which we denote by $J_1$, can be estimated 
using Young's inequality and~\eqref{diff8}. 
In this way, we readily obtain that
\begin{equation}
\label{diff10}
J_1\,\le\,\frac 14\intQt|z^h|^2+\frac 14\intSt|z^h_\Gamma|^2+C\|h\|^2_{H^1(0,t;L^3(\Omega)^3)}\,.
\end{equation}
The remaining two terms, which we denote by $J_2$ and $J_3$, can be handled using the 
H\"older, Young and Sobolev inequalities \pier{as well as} \eqref{stabu}. Indeed, we have that
\begin{align}
\label{diff11}
&J_2\,\le\,C\int_0^t\|z^h(s)\|_2\,\|\rho^h(s)-\rhob(s)\|_4^2\,ds\,\le\,\frac 14\intQt|z^h|^2
+C\pier{\int_0^t \|(\rho^h-\rhob)(s)\|^4_{V} ds}\non\\
&\quad\le \,\frac 14\intQt|z^h|^2+C\,\|h\|^4_{H^1(0,t;L^3(\Omega)^3)}\,.
\end{align}
Similar reasoning yields that
\begin{equation}
\label{diff12}
J_3\,\le\, \frac 14\intSt|z^h_\Gamma|^2+C\,\|h\|^4_{H^1(0,t;L^3(\Omega)^3)}\,.
\end{equation}
Therefore, combining the \pier{estimate \eqref{diff8} with \eqref{diff9}--\eqref{diff12}}, we can conclude that
\begin{equation}
\|(z^h,z^h_\Gamma)\|^2_{L^2(0,t;\calV)}\,\le\,C_3\,\|h\|^4_{H^1(0,t;L^3(\Omega)^3)} \quad
\mbox{for all }\,t\in (0,T]\,.
\end{equation}
In conclusion, the inequality \eqref{Frechet} is fulfilled with the choice 
$\,Z(\lambda):=(C_2+C_3)\lambda^4$. The assertion is thus proved.\QED

\vspace{5mm}
\gian{%
With the differentiability shown, the road is paved to derive a first-order necessary
optimality condition for the control problem under investigation.
Indeed, a standard argument 
(which we do not repeat here) 
invoking the chain rule for \Frechet\ derivatives
and the convexity of the admissible set ${\cal U}_{\rm ad}$ 
yields the result stated below, 
where the following abbreviations are used:
\Bsist
  && \phi_1 := \beta_1 (\mub-\hmuQ) , \quad
  \phi_2 := \beta_2 (\muGb-\hmuS), 
  \label{defunodue}
  \\
  && \phi_3 := \beta_3 (\rhob-\hrhoQ) , \quad
  \phi_4 := \beta_4 (\rhoGb-\hrhoS),
  \label{deftrequattro}
  \\
  && \phi_5 := \beta_5 (\rhob(T)-\hrhoO) , \quad
  \phi_6 := \beta_6 (\rhoGb(T)-\hrhoG) .
  \label{defcinquesei}
\Esist
}%

\gian{%
\Bcor
\label{Badnecessary}
Let the assumptions {\em (A1)--(A5)}  be satisfied, and assume that $\ub\in \Uad$ is a solution
to the control problem {\bf (CP)} with associated state $((\mub,\muGb),(\rhob,\rhoGb))=\calS(\ub)$.
Then, with the notation \accorpa{defunodue}{defcinquesei}, we have~that
\begin{align}
  & \intQ \phi_1 \, \eta
  + \intS \phi_2 \, \etaG
  + \intQ \phi_3 \, \xi
  + \intS \phi_4 \, \xiG
  + \iO \phi_5 \, \xi(T)
  \non
  \\
  & \quad {}
  + \iG \phi_6 \, \xiG(T)
  + \beta_7 \intQ \ub \cdot (v-\ub)
  \geq 0
  \quad \hbox{for every $v\in\Uad$},
  \label{badnecessary}
\end{align}
where, for $v\in\Uad$, $(\xi,\xiG,\eta,\etaG)$
is the solution to the linearized problem corresponding to $h:=v-\ub$.
\Ecor
}%


\section{The optimal control problem}
\label{ADJOINT}
\setcounter{equation}{0}

In this section, we examine deeply the control problem {\bf (CP)} 
of minimizing the functional \eqref{cost}
under the control constraint $u\in\Uad$ and \juerg{the state constraint} \State.
First of all, we show the existence of an optimal control.
Then, we eliminate the solution to the linearized problem 
from the necessary condition~\eqref{badnecessary} already established
(with the notations \accorpa{deftrequattro}{defcinquesei}),
by making use of the solution to a proper adjoint problem.
As for the \juerg{first aim}, we have the following result:

\Bthm
\label{Optimum}
Suppose that the assumptions {\rm (A1)--(A5)} hold true. 
Then the optimal control problem {\bf (CP)} has at least one solution, that is,
there exists \juerg{some} $\ub\in\Uad$ such~that
\Beq
  \calJ(\mub,\muGb,\rhob,\rhoGb,\ub) 
  \leq \calJ(\mu,\muG,\rho,\rhoG,u) 
  \quad \hbox{for every $u\in\Uad$},
  \label{optimum}
\Eeq
where $((\mub,\muGb),(\rhob,\rhoGb))$
and $((\mu,\muG),(\rho,\rhoG))$
are the solutions to the state system \State\ corresponding to the controls
$\ub$ and~$u$, respectively.
\Ethm

\proof
We use the direct method.
Thus, we fix a minimizing sequence,
i.e., a~sequence $\{\un\}$ of admissible controls such~that
\Beq
  \lim_{n\to\infty} \calJ(\mun,\muGn,\rhon,\rhoGn,\un)
  = \Lambda := \inf\calJ(\mu,\muG,\rho,\rhoG,u), 
  \label{minimizing}
\Eeq
where the infimum is taken over the set of \juerg{quintuples} $(\mu,\muG,\rho,\rhoG,u)$
\juerg{that satisfy $\,((\mu,\muG),$ $(\rho,\rhoG))=\calS(u)\,$ for some} $u\in\Uad$.
By Theorem~\ref{WPstate}, the estimates \accorpa{separ}{ssb1} hold true
with constants $\rhomin,\rhomax\in(-1,1)$ and $K_1>0$ that do not depend on~$n$.
On the other hand, every $\un$ belongs to~$\Uad$.
Therefore, we have for a subsequence (still indexed by~$n$)
\Bsist
  & \un \to \ub
  & \quad \hbox{weakly star in $(\LQ\infty\cap\H1{\Lx3})^3$},
  \non
  \\
  & (\mun,\muGn) \to (\mub,\muGb)
  & \quad \hbox{weakly star in $\L\infty\calW$},
  \non
  \\
  & (\rhon,\rhoGn) \to (\rhob,\rhoGb)
  & \quad \hbox{weakly star in $\W{1,\infty}\calH\cap\H1\calV\cap\L\infty\calW$}
  \non
  \\
  &&\quad \hbox{and strongly in $\LQ2\times\LS2$}.
  \non
\Esist
\gianni{By Remark~\ref{Generaltrace}}, we can infer that $\ub\in\Uad$.
Moreover, $(f'(\rhon),\fG'(\rhoGn))$ converges to $(f'(\rhob),\fG'(\rhob))$
strongly in $\LQ2\times\LS2$ and
$\rhon\un$ converges to $\rhob\,\ub$ weakly in~$(\LQ2)^3$.
Hence, it is \sfw\ to verify that $((\mub,\muGb),(\rhob,\rhoGb))$
solves the integrated version of the state system \State\
with $u=\ub$ and time-dependent test functions $(v,\vG)\in\L2\calV$,
that is, the system itself.
Finally, we have, by semicontinuity and~\eqref{minimizing},
\Beq
  \calJ(\mub,\muGb,\rhob,\rhoGb,\ub) 
  \leq \lim_{n\to\infty} \calJ(\mun,\muGn,\rhon,\rhoGn,\un)
  = \Lambda \,.
  \non
\Eeq
Therefore, $\ub$ is an optimal control.\QED

The final step consists in eliminating the solution to the linearized problem 
from the necessary condition~\eqref{badnecessary}, 
with the notations \accorpa{deftrequattro}{defcinquesei},
by using the solution to a proper adjoint problem.
However, we cannot deal with the general case, unfortunately.
Indeed, we are \juerg{only} able to treat a slightly less general situation, namely \juerg{when}
\Beq
  \beta_1 = \beta_2 = 0 
  \label{hpsimplification}
\Eeq
\pier{(cf., e.g., \cite{CGS3} for a similar case).}
Furthermore, for a given optimal control~$\ub$,
if \pier{we let} $((\mub,\muGb),(\rhob,\rhoGb))=\calS(\ub)$ 
\pier{be} the corresponding optimal state,
we still keep the notations \accorpa{defunodue}{defcinquesei}\juerg{, noticing}
 that $\phi_1=0$ and $\phi_2=0$ due to~\eqref{hpsimplification},
and also introduce for brevity
\Beq
  \psi := f''(\rhob)
  \aand 
  \psiG := \fG''(\rhoGb) \,.
  \label{defpsi}
\Eeq
Then the adjoint problem \juerg{reads as follows}:
we look for a quadruplet $(p,\pG,q,\qG)$ satisfying the regularity requirements
\begin{align}
  & (p,\pG) \in \pier{\L\infty\calV} \,, \quad
  (q,\qG) \in \L\infty\calH \cap \L2\calV\,,
  \label{regpq}
  \\
  & (p+\tauO q,\pG+\tauG\qG) \in \H1\calVp\,, 
  \label{regptauq}
\end{align}
and solving 
\begin{align}
  & - \< \dt (p+\tauO q,\pG+\tauG\qG) , (v,\vG) >_{\calV}
  + \iO \nabla q \cdot \nabla v
  + \iG \nablaG\qG \cdot \nablaG\vG
  \non
  \\
  & \quad {}
  + \iO \psi q v
  + \iG \psiG \qG \vG
  - \iO \ub \cdot \nabla p \, v
  = \iO \phi_3 \, v
  + \iG \phi_4 \, \vG
  \non
  \\
  & \quad \hbox{\aet\ and for every $(v,\vG)\in\calV$},
  \label{primaA}
  \\
  \separa
  & \iO \nabla p \cdot \nabla v
  + \iG \nablaG\pG \cdot \nablaG\vG
  = \iO q v 
  + \iG \qG \vG
  \non
  \\
  & \quad \hbox{\aet\ and for every $(v,\vG)\in\calV$},
    \label{secondaA}
  \\
  & \< \gianni{(p+\tauO q,\pG+\tauG\qG)}(T) , (v,\vG) >_{\calV}
  = \iO \phi_5 \, v
  + \iG \phi_6 \vG
  \non
  \\
  & \quad \hbox{for every $(v,\vG)\in\calV$} .
  \label{cauchyA}
\end{align}
\Accorpa\Adjoint primaA cauchyA
We notice that the system \accorpa{primaA}{secondaA}
is the variational formulation of the following boundary value problem:
\begin{align}
  & - \dt(p+\tauO q) - \Delta q + \psi q - \ub \cdot \nabla p = \phi_3
  \aand
  - \Delta p = q
  \quad \hbox{in $Q$},
  \non
  \\[1mm]
  & \juerg{- \dt(\pG+\tauG\qG) + \dn q - \DeltaG\qG + \psiG \qG = \phi_4\,,\quad
  \dn\pG - \DeltaG\pG = \qG \,,}
  \non\\
  &\juerg{\mbox{$\,\,p_{|\Sigma}=\pG$ \,and\, $q_{|\Sigma}=\qG$ }\,\quad \hbox{on $\Sigma$}}.
  \non
\end{align}
However, we only use the weak formulation \Adjoint.

\medskip

We discuss \juerg{the well-posedness of this} problem.
We prepare our existence result by solving an approximating problem
depending on a small parameter $\eps\in(0,1)$.
We recall that $\ub$ belongs to $\Uad$.
However, our results \juerg{are valid} under the weaker assumption
\Beq
  \ub \in \juerg{L^\infty(0,T;L^3(\Omega)^3)} .
  \label{hpuA} 
\Eeq
We replace $\ub$ in \eqref{primaA} by the bounded function $\ueps$ defined \aeQ\
by the conditions
\Beq
  \ueps = \ub
  \quad \hbox{where $|\ub|\leq 1/\eps$}
  \aand
  \ueps = \frac 1\eps \, \frac \ub{|\ub|}
  \quad \hbox{where $|\ub|> 1/\eps\,$}.
  \label{defueps}
\Eeq
Moreover, we introduce a viscosity term in~\eqref{secondaA}.
Finally, we approximate the pair $(\phi_5,\phi_6)\in\calH$ by pairs $(\phieps_5,\phieps_6)$ satisfying
\Beq
  (\phieps_5/\tauO,\phieps_6/\tauG) \in \calV
  \quad \hbox{for $\eps\in(0,1)$}
  \aand
  (\phieps_5,\phieps_6) \to (\phi_5,\phi_6)
  \quad \hbox{in $\calH$ as $\eps\searrow 0$}.
  \quad
  \label{datoeps}
\Eeq
The problem we consider is the following:
we look for a quadruplet $(\peps,\pGeps,\qeps,\qGeps)$
satisfying the regularity requirements
\begin{align}
  (\peps,\pGeps),\, (\qeps,\qGeps) \in \H1\calH \cap \L\infty\calV
  \label{regAeps}
\end{align}
and solving
\begin{align}
  &  - \iO \dt(\peps+\tauO\qeps) v 
  - \iG \dt(\pGeps+\tauG\qGeps) \vG
  + \iO \nabla\qeps \cdot \nabla v
  + \iG \nablaG\qGeps \cdot \nablaG\vG
  \non
  \\
  & \quad {}
  + \iO \psi \qeps v
  + \iG \psiG \qGeps \vG
  - \iO \ueps \cdot \nabla\peps \, v
  = \iO \phi_3 \, v
  + \iG \phi_4 \, \vG\,,
  \label{primaAeps}
  \\
  \separa
  & - \, \eps \iO \dt\peps \, v
  - \eps \iG \dt\pGeps \, \vG
  + \iO \nabla\peps \cdot \nabla v
  + \iG \nablaG\pGeps \cdot \nablaG\vG
  \non
  \\
  & = \iO \qeps v 
  + \iG \qGeps \vG\,,
  \label{secondaAeps}
  \\
  & (\peps,\pGeps)(T) = (0 ,0)
  \aand
  (\qeps,\qGeps)(T) = (\phieps_5/\tauO\,,\phieps_6/\tauG)\,,
  \qquad
  \label{cauchyAeps}
\end{align}
\Accorpa\Adjointeps primaAeps cauchyAeps
where the equalities \accorpa{primaAeps}{secondaAeps} have to hold for every $(v,\vG)\in\calV$ and \aet.
In order to solve this problem, we need a \juerg{preparatory} lemma.

\Blem
\label{Aux}
Let $(\VV,\HH,\VV^*)$ be a Hilbert triplet with $\VV$ separable,
and let the operators 
$\AA\in\calL(\HH,\HH)$, $\BB\in\calL(\VV,\VV^*)$ and $\CC(t)\in\calL(\VV,\HH)$ satisfy,
for some positive constants $\alpha$, $\lambda$ and~$K$,
\begin{align}
  & (\AA w,w)_{\HH}
  \geq \alpha \norma w_{\HH}^2
  \quad \hbox{for every $w\in\HH$}\,,
  \label{coercA}
  \\
  & \< \BB w,w >_{\VV} + \lambda \norma w_{\HH}^2
  \geq \alpha \norma w_{\VV}^2
  \quad \hbox{for every $w\in\VV$}\,,
  \label{coercB}
  \\
  & \norma{\CC(t) w}_{\HH}
  \leq K \norma w_{\VV} 
  \quad \hbox{\aat\ and every $w\in\VV$}\,,
  \label{bddC}
  \\
  & \hbox{for every $w_1\in\VV$ and $w_2\in\HH$},
  \non
  \\
  & \quad \hbox{the function} \quad
  t \mapsto (\CC(t) w_1,w_2)_{\HH}
  \quad \hbox{is measurable on $(0,T)$}.
  \label{measC}
\end{align}
Moreover, assume that $\BB$ is symmetric.
Then, for every $F\in\L2\HH$ and $w_T\in\VV$,
there \juerg{exists} a unique
\Beq
  w \in \H1\HH \cap \L\infty\VV
  \label{regw}
\Eeq
satisfying
\begin{align}
  & - \AA w'(t) + \BB w(t) + \CC(t) w(t) = F(t)
  \quad \hbox{in $\VV^*$ \ \aat}\,,
  \label{astrattaaux}
  \\
  & w(T) = w_T \,.
  \label{finaleaux}
\end{align}
\Elem

\proof
Even nonlinear generalizations of such a result should be known
(see, e.g., \cite{CoVi} for a nonlinear case with $\CC=0$).
However, we did not find any reference that precisely deals with our assumptions.
Therefore, we sketch a short proof.
Both existence and uniqueness are based on the estimate 
obtained by formally testing \eqref{astrattaaux}, written at the time $s$, by $-w'(s)$
and integrating over~$(t,T)$.
By doing this, using the symmetry of~$\BB$
and adding the same quantity to both sides, we obtain 
\Bsist
  && \int_t^T \bigl( \AA w'(s),w'(s) \bigr)_{\HH} \, ds
  + \frac 12 \, \< \BB w(t),w(t) >\juerg{_\VV}
  + \frac \lambda 2 \, \norma{w(t)}_{\HH}^2
  \non
  \\
  && = \frac 12 \, \< \BB w_T,w_T >\juerg{_\VV}
  + \frac \lambda 2 \, \norma{w_T}_{\HH}^2
  + \lambda \int_t^T \bigl( w(s),w'(s) \bigr)_{\HH} \, ds
  \non
  \\
  && \quad {}
  + \int_t^T \bigl( F(s) - \CC(s)w(s) , w'(s) \bigr)_{\HH} \, ds \,.
  \non
\Esist
At this point, we account for \accorpa{coercA}{bddC}, the Young inequality and the Gronwall lemma.
We conclude~that
\Beq
  \norma w_{\H1\HH\cap\L\infty\VV} 
  \leq C \bigl( \norma F_{\L2\HH} + \norma{w_T}_{\VV} \bigr),
  \non
\Eeq
where $C$ depends only on the structural constants and~$T$.
This estimate corresponds to the regularity~\eqref{regw}
and \juerg{implies} that $w=0$ if the data vanish.
However, this is formal, as said at the very beginning.
To make the existence proof rigorous,
we can owe to the separability \juerg{of} $\VV$ and use a Faedo-Galerkin scheme.
To obtain uniqueness, we test \eqref{astrattaaux} by the function $-w'_\delta$ rather than
\juerg{by} ~$-w'$,
where $w_\delta\in\H1\VV$ is obtained by solving the abstract elliptic problem
(here $\II:\VV\to\VV^*$ is the injection)
\Beq
  w_\delta(t) + \delta (\BB+\lambda\II) w_\delta(t) = w(t)
  \quad \aat .
  \non
\Eeq
Then, we use \cite[Appendix: Prop.~6.1-6.3 and Rem.~6.4]{CGG3} in letting $\delta$ tend to zero.
This yields the desired estimate, thus uniqueness if $F=0$ and $w_T=0$.\QED

\Bthm
\label{WellposednessAeps}
For every $\eps\in(0,1)$, the approximating problem \Adjointeps\ has a unique solution 
$(\peps,\pGeps,\qeps,\qGeps)$ satisfying \eqref{regAeps}.
\Ethm

\proof
We present the problem in a different form.
We term $(z,\zG)$ rather than $(v,\vG)$ the arbitrary element of $\calV$ that appears in \eqref{primaAeps}
and add this equation to \eqref{secondaAeps} divided by~$\eps\tau$,
where
\Beq
  \tau := \min \{ \tauO \,, \tauG \}.
  \label{deftau}
\Eeq
This yields \juerg{the identity}
\Bsist
  && - \, \frac 1\tau \iO \dt\peps \, v
  - \frac 1\tau \iG \dt\pGeps \, \vG
  + \frac 1 {\eps\tau} \iO \nabla\peps \cdot \nabla v
  + \frac 1 {\eps\tau} \iG \nablaG\pGeps \cdot \nablaG\vG
  \non
  \\
  && \quad {}
  - \iO \dt(\peps+\tauO\qeps) z 
  - \iG \dt(\pGeps+\tauG\qGeps) \zG
  + \iO \nabla\qeps \cdot \nabla z
  + \iG \nablaG\qGeps \cdot \nablaG\zG
  \non
  \\
  && \quad {}
  + \iO \psi \qeps z
  + \iG \psiG \qGeps \zG
  - \iO \ueps \cdot \nabla\peps \, z
  \non
  \\
  && = \frac 1 {\eps\tau} \iO \qeps v 
  + \frac 1 {\eps\tau} \iG \qGeps \vG
  + \iO \phi_3 \, z
  + \iG \phi_4 \, \zG
  \quad \hbox{for every $(v,\vG),(z,\zG)\in\calV$}.
  \qquad\quad
  \label{astratta}
\Esist
As the pairs $(v,\vG)$ and $(z,\zG)$ are independent from each other,
this equation is equivalent to \accorpa{primaAeps}{secondaAeps},
and we are going to transform it into an abstract equation like \eqref{astrattaaux}
in the framework of the Hilbert triplet
\Beq
  \VV \subset \HH \subset \VVp 
  \quad \hbox{where} \quad
  \VV := \calV \times \calV
  \aand
  \HH := \calH \times \calH \,,
  \quad \hbox{whence} \quad
  \VVp = \calV^*\times\calV^* 
  \non
\Eeq
with a non-standard embedding $\HH\subset\VVp$,
due to a particular choice of the inner product in~$\HH$.
In order to simplify the notation, we write the elements $((v,\vG),(z,\zG))$ of $\HH$
as quadruplets $(v,\vG,z,\zG)$.
We~set
\Bsist
  && \bigl( (p,\pG,q,\qG) , (v,\vG,z,\zG) \bigr)_{\HH}
  := \iO \Bigl( \frac 1 {\eps\tauO} \, pv + qz \Bigr)
  + \iG \Bigl( \frac 1 {\eps\tauG} \, \pG\vG + \qG\zG \Bigr)
  \non
  \\
  && \quad \hbox{for every $(p,\pG,q,\qG),\,(v,\vG,z,\zG)\in\HH$}
  \label{prodH}
\Esist
and notice that $(\cpto,\cpto)_{\HH}$ actually is an inner product 
and that the corresponding norm is equivalent to the standard one.
Moreover, we define the operators
$\AA^\eps\in\calL(\HH,\HH)$, $\BB\in\calL(\VV,\VV^*)$ and $\CC^\eps(t)\in\calL(\VV,\HH)$ 
by the formulas
\Bsist
  && \AA^\eps (p,\pG,q,\qG)
  := \Bigl( \frac {\eps\tauO}\tau \, p , \frac {\eps\tauG}\tau \, \pG , p+\tauO q , \pG+\tauG\qG \Bigr)
  \non
  \\
  && \quad \hbox{for every $(p,\pG,q,\qG)\in\HH$}\,,
  \label{defAeps}
  \\[1mm]
  && \< \gianni{\BB^\eps} (p,\pG,q,\qG) , (v,\vG,z,\zG) >_{\VV}
  \non
  \\
  && \,\,:= \frac 1 {\eps\tau} \iO \nabla p \cdot \nabla v
  + \frac 1 {\eps\tau} \iG \nablaG\pG \cdot \nablaG\vG
  + \iO \nabla q \cdot \nabla z
  + \iG \nablaG\qG \cdot \nablaG\zG
  \non
  \\
  && \quad \hbox{for every $(p,\pG,q,\qG),\,(v,\vG,z,\zG)\in\VV$}\,,
  \label{defB}
  \\[1mm]
  && \CC^\eps(t) (p,\pG,q,\qG)
  \non
  \\
  && \,\,:= \Bigl( 
    - \frac \tauO\tau \, q \,,\,
    - \frac \tauG\tau \, \qG \,,\,
    -\ueps(t) \cdot \nabla p + \psi(t) \, q \,,\,
    \psiG(t) \, \qG
  \Bigr)
  \non
  \\
  && \quad \hbox{\aat\ and every $(p,\pG,q,\qG)\in\VV$} \,.
  \label{defCeps}
\Esist
A simple computation shows that
\Bsist
  && \bigl( \AA^\eps (p,\pG,q,\qG) , (v,\vG,z,\zG) \bigr)_{\HH}
  \non
  \\
  && \,\,= \iO \Bigl( \frac 1\tau \, pv + (p+\tauO q) z \Bigr)
  + \iG \Bigl( \frac 1\tau \, \pG\vG + (\pG+\tauG\qG) \pier{\zG} \Bigr)
  \non
  \\
  && \quad \hbox{for every $(p,\pG,q,\qG),\,(v,\vG,z,\zG)\in\HH$}\,,
  \non
  \\
  && \bigl( \CC^\eps(t) (p,\pG,q,\qG) , (v,\vG,z,\zG) \bigr)_{\HH}
  \non
  \\
  && \,\,= - \frac 1 {\eps\tau} \iO q v 
  - \frac 1 {\eps\tau} \iG \qG \vG
  + \iO \psi(t) q z
  + \iG \psiG(t) \qG \zG
  - \iO \ueps(t) \cdot \nabla p \, z
  \non
  \\
  && \quad \hbox{\aat\ and every $(p,\pG,q,\qG)\in\VV$ and $(v,\vG,z,\zG)\in\HH$}.
  \non
\Esist
Therefore, the variational equation \eqref{astratta} takes the form
\Bsist
  && - \bigl( \AA^\eps \dt(\peps,\pGeps,\qeps,\qGeps)(t) , (v,\vG,z,\zG) \bigr)_{\HH}
  + \< \gianni{\BB^\eps} (\peps,\pGeps,\qeps,\qGeps)(t) , (v,\vG,z,\zG) >_{\VV}
  \non
  \\
  && \quad {}
  + \bigl( \CC^\eps(t) (\peps,\pGeps,\qeps,\qGeps)(t) , (v,\vG,z,\zG) \bigr)_{\HH}
  = \bigl( F(t) , (v,\vG,z,\zG) \bigr)_{\HH}
  \non
  \\
  && \quad \hbox{\aat\ and every $(v,\vG,z,\zG)\in\VV$}\,,  \non
\Esist
with an obvious definition of $F\in\L2\HH$.
Thus, it is a particular case of~\eqref{astrattaaux}.
On the other hand, \eqref{cauchyAeps} is equivalent~to
\Bsist
  && \bigl( (\peps,\pGeps,\qeps,\qGeps)(T) , (v,\vG,z,\zG) \bigr)_{\HH}
  = \bigl( (0,0,\phieps_5/\tauO,\phieps_6/\tauG) , (v,\vG,z,\zG) \bigr)_{\HH}
  \non
  \\
  && \quad \hbox{for every $(v,\vG,z,\zG)\in\HH$}\,,
  \non
\Esist
and $(0,0,\phieps_5/\tauO,\phieps_6/\tauG)$ belongs to~$\VV$ by the 
first \pier{condition in}~\eqref{datoeps}.
Therefore, in order to conclude, it is sufficient to check the properties \accorpa{coercA}{measC}
of Lemma~\ref{Aux}.
The second and fourth ones are clear, and $\gianni{\BB^\eps}$ is obviously symmetric.
Moreover, \eqref{bddC}~easily follows from the boundedness of $\ueps$, $\psi$ and~$\psiG$.
As for~\eqref{coercA}, we have, for some constant $\alpha>0$ and every $(v,\vG,z,\zG)\in\HH$,
\Bsist
  && \bigl( \AA^\eps (v,\vG,z,\zG) , (v,\vG,z,\zG) \bigr)_{\HH}
  \non
  \\
  && = \iO \Bigl( \frac 1\tau \, |v|^2 + (v+\tauO z)z \Bigr)
  + \iG \Bigl( \frac 1\tau \, |\vG|^2 + (\vG+\tauG\zG)\zG \Bigr)
  \non
  \\
  && \geq \iO \Bigl( \frac 1\tau \, |v|^2 + \tau |z|^2 + vz \Bigr)
  + \iG \Bigl( \frac 1\tau \, |\vG|^2 + \tau |\zG|^2 + \vG\zG \Bigr)
  \non
  \\
  && \geq \frac 12 \iO \Bigl( \frac 1\tau \, |v|^2 + \tau |z|^2 \Bigr)
  + \frac 12 \iG \Bigl( \frac 1\tau \, |\vG|^2 + \tau |\zG|^2 \Bigr)
  \geq \alpha \norma{(v,\vG,z,\zG)}_{\HH}^2\,,
  \non
\Esist
the last inequalities being due to the Young inequality 
and to the equivalence between the norm $\norma\cpto_{\HH}$ induced by the inner product \eqref{prodH}
and the natural norm of~$\HH$.
Therefore, Lemma~\ref{Aux} can be applied and the proof is complete.\QED

\Bthm
\label{WellposednessA}
Let the assumptions {\em (A1)--(A5)} and \eqref{hpsimplification} be satisfied. 
Moreover, assume that $\ub\in \Uad$ is a solution to the control problem {\bf (CP)} 
and that $((\mub,\muGb),(\rhob,\rhoGb))=\calS(\ub)$ is the associated state. 
Then, with the notations \accorpa{deftrequattro}{defcinquesei} and~\eqref{defpsi}, 
the adjoint problem \Adjoint\ has \pier{a unique} solution $(p,\pG,q,\qG)$
satisfying \accorpa{regpq}{regptauq}.
\Ethm

\proof
In order to show existence, 
we perform a number of a priori estimates on the solution 
$(\peps,\pGeps,\qeps,\qGeps)$ to the approximating problem.
However, we explicitly write the superscript $\eps$ only at the end of each estimate.
Moreover, we make use of the same symbol $c$ 
to denote different constants that do not depend on~$\eps$.
The symbol $c_\delta$ \juerg{stands} for (possibly different) constants 
that can also depend on the parameter~$\delta$.
By denoting by $|\Omega|$ and $|\Gamma|$ the volume of $\Omega$ and the area of~$\Gamma$, respectively, 
we define \juerg{the} mean value function~by
\Beq
  \mean(v,\vG)
  := \frac {\iO v + \iG \vG} {|\Omega|+|\Gamma|} 
  \quad \hbox{for every $(v,\vG)\in\calH$}
  \label{defmean}
\Eeq
and observe that the Poincar\'e type inequality 
\Beq
  \normaVV{(v,\vG)}
  \leq c \bigl(
    \norma{(\nabla v,\nablaG\vG)}_{\pier{\calH^3}}
    + |\mean(v,\vG)|
  \bigr)
  \quad \hbox{for every $(v,\vG)\in\calV$}
  \label{poincare}
\Eeq
holds true with a constant $c$ that depends only on~$\Omega$.
Finally, we set, for brevity,
\Beq
  Q^t:=\Omega\times(t,T)
  \aand
  \Sigma^t:=\Gamma\times(t,T)
  \quad \hbox{for $t\in(0,T)$}.
  \non
\Eeq

\step
First a priori estimate

We test \eqref{secondaAeps} by $(1,1)$ and obtain
\Beq
  |\mean(q,\qG)| \, (|\Omega|+|\Gamma|)
  \leq \eps  \iO |\dt p| + \eps \iG |\dt\pier{\pG}| 
  \quad \aet.
\Eeq
As $\eps^2\leq\eps$, since $\eps\in(0,1)$, we infer that
\Beq
  \int_t^T |\mean(q,\qG)(s)|^2 \, ds
  \leq c\,\eps \bintQt |\dt p|^2 + c\,\eps \bintSt |\dt\pier{\pG}|^2 
  \quad \hbox{for every $t\in[0,T]$}
  \non
\Eeq
where $c$ depends only on~$\Omega$.
On the other hand, \pier{owing also} to \eqref{poincare} and to the Sobolev inequality, we deduce that
\Bsist
  && \int_t^T \norma{q(s)}_6^2 \, ds
  \leq c \Bigl(
    \bintQt |q|^2
    + \bintQt |\nabla q|^2
  \Bigr)
  \non
  \\
  && \leq c \Bigl(
    \bintQt |\nabla q|^2
    + \bintSt |\nablaG\qG|^2
    + \int_t^T |\mean(q,\qG)(s)|^2 \, ds
  \Bigr) .
  \non
\Esist
Therefore, by combining these inequalities, we conclude that
\Beq
  \int_t^T \norma{q(s)}_6^2 \, ds
  \leq \CO \Bigl(
    \bintQt |\nabla q|^2
    + \bintSt |\nablaG\qG|^2  
    + \eps \bintQt |\dt p|^2 + \eps \bintSt |\dt\pier{\pG}|^2
  \Bigr)
  \label{primastimaA}
\Eeq
for every $t\in[0,T]$, with a constant $\CO$ that depends only on~$\Omega$.

\step
Second a priori estimate

We test \eqref{primaAeps} by $(q,\qG)$, integrate over~$(t,T)$,
account for the Cauchy conditions \eqref{cauchyAeps}, and have
\Bsist
  && - \bintQt \dt p \, q
  - \bintSt \dt\pG \, \qG
  + \frac \tauO 2 \iO |q(t)|^2
  + \frac \tauG 2 \iG |\qG(t)|^2
  + \bintQt |\nabla q|^2
  + \bintSt |\nablaG\qG|^2
  \non
  \\
  && = \frac \tauO 2 \iO |\phieps_5/\tauO|^2
  + \frac \tauG 2 \iG |\phieps_6/\tauG|^2 
  - \bintQt \psi |q|^2
  - \bintSt \psiG |\qG|^2
  \non
  \\
  && \quad {}
  + \bintQt \ueps \cdot \nabla p \, q
  + \bintQt \phi_3 \, q
  + \bintSt \phi_4 \, \qG \,.
  \non
\Esist
At the same time, we test \eqref{secondaAeps} by $-\dt(p,\pG)$ \juerg{and} integrate over~$(t,T)$ 
\juerg{to obtain}
\Beq
  \eps \bintQt |\dt p|^2
  + \eps \bintSt |\dt\pG|^2
  + \frac 12 \iO |\nabla p(t)|^2
  + \frac 12 \iG |\nablaG\pG(t)|^2
  = - \bintQt q \dt p
  - \bintSt \qG \dt\pG \,.
  \non
\Eeq
Now, we add this equality to the previous one 
and observe that four terms cancel each other.
Moreover, \pier{accounting for \eqref{defueps} and \eqref{primastimaA}, we treat the transport term as follows:} 
\Bsist
  && \bintQt \ueps \cdot \nabla p \, q
  \leq \norma\ub_\juerg{{L^\infty(0,T;L^3(\Omega)^3)}} \int_t^T \norma{\nabla p(s)}_2 \, \norma{q(s)}_6 \, ds
  \non
  \\
  && \leq \delta \int_t^T \norma{q(s)}_6^2 \, ds 
  + c_\delta \bintQt |\nabla p|^2 
  \non
  \\
  && \leq \delta \, \CO \Bigl(
    \bintQt |\nabla q|^2
    + \bintSt |\nablaG\qG|^2  
    + \eps \bintQt |\dt p|^2 + \eps \bintSt |\dt\pier{\pG}|^2
  \Bigr)
  + c_\delta \bintQt |\nabla p|^2\,,
  \non
\Esist
\pier{where $\delta>0$ is arbitrary.} 
Therefore, since $\psi$ and~$\psiG$ are bounded, $\phi_3$ and $\phi_4$ are $L^2$-functions,
and \eqref{datoeps} implies that $(\phieps_5,\phieps_6)$ is bounded in $\calH$
uniformly with respect to~$\eps$,
by choosing $\delta$ such that $\delta\CO\leq1/2$ and using the Gronwall lemma,
we conclude that
\Beq
  \norma{(\qeps,\qGeps)}_{\L\infty\calH\cap\L2\calV} 
  + \norma{(\nabla\peps,\nablaG\pGeps)}_{\pier{\L\infty\calH}}
  + \eps^{1/2} \norma{\dt(\peps,\pGeps)}_{\L2\calH}
  \leq c \,.
  \label{secondastimaA}
\Eeq

\step
Third a priori estimate

By testing \eqref{primaAeps} by an arbitrary pair $(v,\vG)\in\L2\calV$
and accounting for~\eqref{secondastimaA},
we easily deduce that
\Beq
  \norma{\dt(\peps+\tauO\qeps,\pGeps+\tauG\qGeps)}_{\L2\calVp}
  \leq c \,.
  \label{terzastimaA}
\Eeq

\step
Fourth a priori estimate

Clearly, \eqref{terzastimaA} implies that
\Beq
  \norma{{\textstyle\frac d{dt}}\mean(p+\tauO q,\pG+\tauG\qG)}_{L^2(0,T)}
  \leq c \norma{\dt(p+\tauO q,\pG+\tauG\qG)}_{\L2\calVp} \, \normaVV{(1,1)}
  \leq c \,,
  \non
\Eeq
and\pier{, in view of \eqref{cauchyAeps},} we infer that
\Beq
  \norma{\mean(p+\tauO q,\pG+\tauG\qG)}_{L^\infty(0,T)}
  \leq c \,.
  \label{perquartastimaA}
\Eeq
On the other hand, even $(\tauO q,\tauG\qG)$ is bounded in $\L\infty\calH$ by~\eqref{secondastimaA},
and \pier{consequently} $\mean(\tauO q,\tauG\qG)$ is bounded in $L^\infty(0,T)$. 
Therefore, \eqref{perquartastimaA} ensures that the same holds for $\mean(p,\pG)$.
By accounting for \pier{\eqref{secondastimaA}} and the Poincar\'e type inequality~\eqref{poincare}, 
we conclude that
\Beq
  \norma{(\peps,\pGeps)}_{\pier{\L\infty\calV}} \leq c \,.
  \label{quartastimaA}
\Eeq

\step
Existence

We are ready to take the limit as $\eps\searrow 0$.
We have, at least for a subsequence,
\begin{align}
  & (\peps,\pGeps) \to (p,\pG)
  \quad \hbox{weakly \pier{star} in $\pier{\L\infty\calV}$}\,,
  \label{wconvp}
  \\
  & (\qeps,\qGeps) \to (q,\qG)
  \quad \hbox{weakly star in $\L\infty\calH\cap\L2\calV$}\,,
  \label{wconvq}
  \\
  & \dt(\peps+\tauO\qeps,\pGeps+\tauG\qGeps) \to \dt(p+\tauO q,\pG+\tauG\qG)
  \quad \hbox{weakly in $\L2\calVp$}\,,
  \label{wconvdt}
  \\
  & \eps \, \dt(\peps,\pGeps) \to 0
  \quad \hbox{strongly in $\L2\calH$} .
  \label{wconvepsdt}
\end{align}
As \accorpa{wconvp}{wconvdt} imply that \pier{(see, e.g., \cite[Sect.~8,~Cor.~4]{Simon})}
\begin{align*}
  &(\peps+\tauO\qeps,\pGeps+\tauG\qGeps) \to (p+\tauO q,\pG+\tauG\qG) \\
  &\quad  \pier{\hbox{strongly in~$\C0\calVp\cap \L2\calH$}}, \non
\end{align*}
and \juerg{since} the approximating final data satisfy the \pier{convergence property in}~\eqref{datoeps}, \pier{from \eqref{cauchyAeps} it follows that}
the Cauchy condition \eqref{cauchyA} is fulfilled by the limiting quadruplet $(p,\pG,q,\qG)$.
Moreover, by recalling~\eqref{defueps},
we see that $\ueps$ converges to $\ub\,$ \aeQ.
By combining this with \eqref{hpuA} and the inequality $|\ueps|\leq|\ub|$ \aeQ, 
we deduce that $\ueps$ converges to $\ub$ strongly (e.g.) in~$(\LQ{8/3})^3$.
Thus, by also accounting for~\eqref{wconvp}, we infer that
$\ueps\cdot\nabla\peps$ converges to $\ub\cdot\nabla p$ weakly in $\LQ{8/7}$.
Therefore, we can take the limit in the integrated version of \accorpa{primaAeps}{secondaAeps} 
with time-dependent test functions $(v,\vG)\in\L2\calV$ with $v\in\LQ8$
and conclude that $(p,\pG,q,\qG)$ solves
the integrated version of \accorpa{primaA}{secondaA} with the same test functions.
By density, since $\ub\cdot\nabla p$ belongs to $\L2{\Lx{6/5}}$ by \eqref{hpuA} and~\eqref{wconvp},
and $\L2V\subset\L2{\Lx6}$ by the Sobolev inequality,
one can take any element of $\L2\calV$ as a test function.
Hence, this quadruplet also solves the equations \eqref{primaA} and \eqref{secondaA} as they~are.

\step  
Uniqueness

\pier{Since the problem \Adjoint\ is linear, it is sufficient to prove that the unique solution
with $(\phi_3,\phi_4,\phi_5,\phi_6)=(0,0,0,0)$ is $(p,\pG,q,\qG)=(0,0,0,0)$.
In the next lines, $C_i$,~$i=1,2,\dots$,
denote positive constants that depend only on the structural assumptions 
and the $L^\infty$-norms of $\ub$, $\psi$ and~$\psiG$.}

\pier{First, we introduce the primitive functions
\begin{align*} 
Q(t) := - \inttT q(s)\, ds , \quad \QG (t) := - \inttT \qG (s)\, ds, \quad \  t\in [0,T],
\end{align*} 
and integrate \eqref{primaA} from $t$ to $T$ in order to obtain
\begin{align}
  &  \< \, (p+\tauO q,\pG+\tauG\qG)(t) , (v,\vG) >_{\calV}
  - \iO \nabla Q (t)\cdot \nabla v
  - \iG \nablaG\QG (t)\cdot \nablaG\vG
  \non
  \\
  & {}
  + \iO \inttT (\psi q)(s) \,ds\, v
  + \iG \inttT (\psiG \qG) (s) \,ds\, \vG
  - \iO \inttT (\ub \cdot \nabla p) (s)\, ds \, v 
  \non
  \\[0.2cm]
  &  = 0 \quad \hbox{ for every $t\in [0,T]$ and $(v,\vG)\in\calV$}.
  \label{pier3}
\end{align}
Next, we test \eqref{secondaA}, written at the time~$t$, by $(p,\pG)(t)$;
at the same time, we take $(v,\vG) = ( q,\qG)(t)$ in \eqref{pier3} and sum the two equalities we obtain  
by observing that there is a cancellation of four terms. Integrating once more with respect to $t$, we deduce that 
\begin{align}
  &\tauO \bintQt |q|^2
  + \tauG \bintSt |\qG|^2
  + \iO |\nabla Q (t)|^2
  \non
  \\  &+ \iG |\nablaG\QG (t)|^2
  + \bintQt |\nabla p|^2
  + \bintSt |\nablaG\pG|^2
  \non
  \\
  &= 
    - \bintQt  \intsT (\psi q)(\sigma) \,d\sigma\, q(s)\,ds
    - \bintSt  \intsT (\psiG \qG) (\sigma) \,d\sigma\, \qG (s)\, ds 
 \non
  \\
    &\quad + \bintQt  \intsT (\ub \cdot \nabla p) (\sigma)\,d\sigma \, q(s)\, ds\,.     
\label{pier4}
\end{align} 
\juerg{We} now estimate the terms on the \rhs\ of \eqref{pier4}.
Thanks to the Young and H\"older inequalities, and using the $L^\infty$-bound for $\psi$, we infer that 
\begin{align}
&- \bintQt  \intsT (\psi q)(\sigma) \,d\sigma\, q(s)\,ds \leq 
\frac{\tauO}4 \bintQt |q|^2 + 
\frac1{\tauO} \norma{\psi}^2_\infty \, T \bintQt \intsT | q(\sigma)|^2 \,d\sigma\, ds
\non\\
&\qquad \leq \frac{\tauO}4 \bintQt |q|^2 + 
C_1 \inttT \bigg( \bintQs | q|^2 \bigg)ds .
\label{pier5}
\end{align}
Arguing similarly for the boundary integral, we have that
\begin{align}
- \bintSt  \intsT (\psiG \qG) (\sigma) \,d\sigma\, \qG (s) \,ds  \leq 
\frac{\tauG}4 \bintSt |\qG|^2 + 
C_2 \inttT \bigg( \bintSs | \qG |^2  \bigg) ds.
\label{pier6}
\end{align}
Also the last term of \eqref{pier4} can be treated by the same token. Indeed, we see that
\begin{align}
&\bintQt  \intsT (\ub \cdot \nabla p) (\sigma) \,d\sigma \, q(s)\, ds  \leq 
\frac{\tauO}4 \bintQt |q|^2 + 
\frac1{\tauO} \norma{\ub}^2_\infty \, T \bintQt \intsT | \nabla p (\sigma)|^2 \,d\sigma\, ds
\non\\
&\qquad \leq \frac{\tauO}4 \bintQt |q|^2 + 
C_3 \inttT \bigg( \bintQs |\nabla p |^2 \bigg)ds .
\label{pier7}
\end{align}
At this point, if we combine the inequality \eqref{pier4} with the 
estimates~\eqref{pier5}--\eqref{pier7} and then apply the Gronwall lemma, we deduce that 
$(q,\qG)=(0,0) $ as well as that the vectors $\nabla p$ and $\nablaG \pG$ vanish. Hence, by a comparison in
\eqref{pier3}, we finally conclude that $(p,\pG)=(0,0)$, and the proof is complete.}\QED

\medskip

Once the solvability of problem \Adjoint\ is established,
we actually can eliminate the solution to the linearized problem in~\eqref{badnecessary},
as stated in our final result.
For its proof, we need a Leibniz rule
which is well known under slightly different assumptions.

\Blem
\label{Leibniz}
Assume that
\Beq
  y \in \H1\calH \cap \L2\calV
  \aand
  z \in \H1\calVp \cap \L2\calH \,.
  \label{hpleibniz}
\Eeq
Then the function $t\mapsto(y(t),z(t))_{\calH}$
is absolutely continuous on~$[0,T]$, 
and its derivative is given~by
\Beq
  \frac d{dt} \, (y,z)_{\calH}
  = (y',z)_{\calH}
  + \< z',y>_{\calV}
  \quad \aet.
  \label{leibniz}
\Eeq
\Elem

\proof
By the trace method with $p=2$ of the interpolation theory 
(see, e.g., \cite{Lun}), the continuous embeddings
\Bsist
  && \H1\calH \cap \L2\calV \subset \C0{(\calV,\calH)_{1/2}}\,,
  \non
  \\
  && \H1\calVp \cap \L2\calH \subset \C0{(\calH,\calVp)_{1/2}}\,,
  \non
\Esist
hold true, as well as the duality formula $(\calH,\calVp)_{1/2}=(\calV,\calH)_{1/2}^*\,$.
Therefore, the map
\Beq
  t \mapsto \bigl( y(t),z(t) \bigr)_{\calH}
  = \< z(t) , y(t) >_{(\calV,\calH)_{1/2}}
  \non
\Eeq
is continuous on~$[0,T]$.
Thus, to conclude, it suffices to prove that
\Bsist
  && \< z(t_2) , y(t_2) >_{(\calV,\calH)_{1/2}}
  - \< z(t_1) , y(t_1) >_{(\calV,\calH)_{1/2}}
  \non
  \\
  && = \int_{t_1}^{t_2} \bigl(
    (y'(s),z(s))_{\calH}
    + \< z'(s),y(s)>_{\calV}
  \bigr) \, ds
  \quad \hbox{for every $t_1,t_2\in[0,T]$}.
  \label{tesileibniz}
\Esist
To this end, we approximate $z$ by functions $z_n\in\H1H$ satisfying
\Beq
  z_n \to z
  \quad \hbox{in $\H1\calVp\cap\L2\calH$}.
  \non
\Eeq
Then, \eqref{tesileibniz} holds for $y$ and~$z_n$, as is well known.
At this point, one lets $n$ tend to infinity and obtains \eqref{tesileibniz}
by observing that $z_n$ converges to $z$ also in $\C0{(\calV,\calH)_{1/2}^*}$.\QED

\Bthm
\label{Goodnecessary}
Let the assumptions {\em (A1)--(A5)} and \eqref{hpsimplification} be satisfied. 
Moreover, assume that $\ub\in \Uad$ is a solution to the control problem {\bf (CP)} 
with associated state $((\mub,\muGb),(\rhob,\rhoGb))=\calS(\ub)$. 
Furthermore, with the notations \accorpa{deftrequattro}{defcinquesei} and~\eqref{defpsi}, 
let $(p,\pG,q,\qG)$ be \pier{the solution} to the adjoint problem \Adjoint\
satisfying the regularity requirements \accorpa{regpq}{regptauq}.
Then, we have  
\Beq
  \intQ \bigl( \rhob \, \nabla p + \beta_7 \ub \bigr) \cdot (v-\ub)
  \geq 0
  \quad \hbox{for every $v\in\Uad\,$}.
  \label{goodnecessary}
\Eeq
\Ethm

\proof
We fix any $v\in\Uad$ and introduce the linearized problem corresponding to the choice $h=v-\ub$
as in Corollary~\ref{Badnecessary}.
Then, we test \eqref{primaL} and \eqref{secondaL} by $(p,\pG)$ and $(q,\qG)$, respectively,
integrate over~$(0,T)$ and sum up.
We obtain
\Bsist
  && \intQ \dt\xi \, p
  + \intS \dt\xiG \, \pG
  + \intQ \nabla\eta \cdot \nabla p
  + \intS \nablaG\etaG \cdot \nablaG\pG
  - \intQ \xi \ub \cdot \nabla p
  - \intS \rhob \, h \cdot \nabla p
  \qquad
  \non
  \\
  \separa
  && + \, \tauO \intQ \dt\xi \, q
  + \tauG \intS \dt\xiG \, \qG
  + \intQ \nabla\xi \cdot \nabla q
  + \intS \nablaG\xiG \cdot \nablaG\qG
  \non
  \\
  && \quad {}
  + \intQ \psi \xi q
  + \intS \psiG \xiG \qG
  = \intQ \eta q 
  + \intS \etaG \qG\,, 
  \label{testedL}
\Esist
and we observe that the sum of the terms \juerg{involving} time derivatives can be written~as
\Bsist
  && \intQ \dt\xi \, p
  + \intS \dt\xiG \, \pG
  + \tauO \intQ \dt\xi \, q
  + \tauG \intS \dt\xiG \, \qG
  \non
  \\
  && = \ioT \bigl( \dt(\xi,\xiG)(s) , (p+\tauO q,\pG+\tauG\qG(s) \bigr)_{\calH} \, ds.
  \non
\Esist
Now, we test \eqref{primaA} and \eqref{secondaA} by $-(\xi,\xiG)$ and $-(\eta,\etaG)$, respectively,
integrate over~$(0,T)$, and sum up.
We obtain \juerg{the identity}
\Bsist
  && \ioT \< \dt(p+\tauO q,\pG+\tauG\qG)(s) , (\xi,\xiG)(s) >_{\calV} \, ds
  - \intQ \nabla q \cdot \nabla\xi
  - \intS \nablaG\qG \cdot \nablaG\xiG
  \non
  \\
  \separa
  && \quad {}
  - \intQ \psi q \xi
  - \intS \psiG \qG \xiG
  + \intQ \ub \cdot \nabla p \, \xi
  \non
  \\
  && - \, \intQ \nabla p \cdot \nabla\eta
  - \intQ \nablaG\pG \cdot \nablaG\etaG
  \non
  \\
  && = - \intQ \phi_3 \, \xi
  - \intS \phi_4 \, \xiG
  - \intQ q \eta 
  - \intS \qG \etaG \,.
  \label{testedA}
\Esist
At this point, we add \juerg{the} equalities \eqref{testedL} and \eqref{testedA} to each other.
Then, the most part of the terms cancels out,
and the sum of the integrals involving time derivatives
can be treated by \pier{invoking} Lemma~\ref{Leibniz}.
Hence, we obtain
\Bsist
  && \ioT \frac d{dt} \, \bigl( (p+\tauO q,\pG+\tauG\qG)(s) , (\xi,\xiG)(s) \bigr)_{\calH} \, ds
  - \intS \rhob \, h \cdot \nabla p
  \non
  \\
  && =
  - \intQ \phi_3 \, \xi
  - \intS \phi_4 \, \xiG \,.
  \non
\Esist
On the other hand, thanks to Lemma~\ref{Leibniz}, \eqref{cauchyL} and~\eqref{cauchyA} with $(v,\vG)=(\xi,\xiG)(T)$, 
we also have
\Bsist
  && \ioT \frac d{dt} \, \bigl( (p+\tauO q,\pG+\tauG\qG)(s) , (\xi,\xiG)(s) \bigr)_{\calH} \, ds
  \non
  \\
  && = \bigl( (p+\tauO q,\pG+\tauG\qG)(T) , (\xi,\xiG)(T) \bigr)_{\calH}
  - \bigl( (p+\tauO q,\pG+\tauG\qG)(0) , (\xi,\xiG)(0) \bigr)_{\calH}
  \non
  \\
  && = \bigl( (\phi_5,\phi_6) , (\xi,\xiG)(T) \bigr)_{\calH} \,.
  \non
\Esist
Therefore, \eqref{badnecessary} becomes~\eqref{goodnecessary}.


\vspace{3truemm}

{\small%

\Begin{thebibliography}{10}

\bibitem{CFGS1}
P. Colli, \pier{M.H. Farshbaf-Shaker, G. Gilardi,} J. Sprekels: 
{\em Optimal boundary control of a 
viscous Cahn-Hilliard system with dynamic boundary condition 
and double obstacle potentials.}
SIAM J. Control Optim. {\bf 53} (2015), 2696-2721.

\bibitem{CFGS2}
P. Colli, M.H. Farshbaf-Shaker, G. Gilardi, J. Sprekels: {\em Second-order analysis of a boundary control problem for the viscous Cahn-Hilliard equation with dynamic boundary conditions.}
\pier{Ann. Acad. Rom. Sci. Ser. Math. Appl.} {\bf 7} (2015), 41-66.

\pier{%
\bibitem{CF1} 
P.\ {C}olli, T.\ {F}ukao: 
{\em {C}ahn--{H}illiard equation with dynamic boundary conditions 
and mass constraint on the boundary.} 
J. Math. Anal. Appl. {\bf 429} (2015), 1190-1213.
\bibitem{CF2} 
P.\ {C}olli, T.\ {F}ukao: 
{\em Equation and dynamic boundary condition of 
Cahn--Hilliard type with singular potentials.}
Nonlinear Anal. {\bf 127} (2015), 413-433.%
}

\bibitem{CGG3}
P. Colli, G. Gilardi, M. Grasselli:
{\em Well-posedness of the weak formulation for the phase-field model with memory.}
\pier{Adv.} Differential Equations {\bf 2} (1997), 487-508.

\bibitem{CGRS1}
P. Colli, G. Gilardi, E. Rocca, J. Sprekels: 
{\em Optimal distributed control of a diffuse interface model of tumor growth.} 
Nonlinearity \pier{{\bf 30} (2017), 2518-2546.}

\bibitem{CGS1}
P. Colli, G. Gilardi, J. Sprekels: {\em A boundary control problem 
for the pure Cahn--Hilliard equation 
with dynamic boundary conditions.}
Adv. Nonlinear Anal. {\bf 4} (2015), 311-325.

 \bibitem{CGS2}
P. Colli, G. Gilardi, J. Sprekels: {\em Distributed optimal control 
of a nonstandard nonlocal phase field system.}
AIMS \pier{Mathematics} {\bf 1} (2016), 246-281.

\bibitem{CGS3}
P. Colli, G. Gilardi, J. Sprekels: {\em A boundary control problem for the viscous Cahn--Hilliard equation with 
dynamic boundary conditions.}
Appl. Math. Optim. {\bf 72} (2016), 195-225.

\pier{\bibitem{CGS3bis}
P. Colli, G. Gilardi, J. Sprekels: {\em Recent results 
on the Cahn--Hilliard equation 
with dynamic boundary conditions.} Vestn. Yuzhno-Ural. 
Gos. Univ., Ser. Mat. Model. Program. {\bf 10} (2017), 5-21.}

\bibitem{CGS4}
P. Colli, G. Gilardi, J. Sprekels: {\em Distributed optimal control 
of a nonstandard nonlocal phase field 
system with double obstacle potential.}
\juerg{Evol. Equ. Control Theory {\bf 6} (2017), 35-58}.

\bibitem{CGS13}
P. Colli, G. Gilardi, J. Sprekels: {\em On a Cahn--Hilliard system with convection
and dynamic boundary conditions}. \pier{Preprint arXiv:1704.05337 [math.AP] (2017), 
pp.~1-34.}

\bibitem{CoVi}
P. Colli, A. Visintin:
{\em On a class of doubly nonlinear evolution equations.} 
Comm. Partial Differential Equations {\bf 15} (1990), 737-756.

\bibitem{FRS}
S. Frigeri, E. Rocca, J. Sprekels: {\em Optimal distributed control 
of a nonlocal Cahn--Hilliard/Navier--Stokes 
system in two dimensions.}
SIAM J. Control. Optim. {\bf 54} (2016), 221-250. 

\pier{%
\bibitem{FY} 
T. Fukao, N. Yamazaki:
{\em A boundary control problem for the equation 
and dynamic boundary condition of Cahn--Hilliard type.} 
To appear in ``Solvability, Regularity, Optimal Control 
of Boundary Value Problems for PDEs'', 
P.~Colli, A.~Favini, E.~Rocca, G.~Schimperna, J.~Sprekels~(eds.), 
Springer INdAM Series.}

\bibitem{HHKW}
M. Hinterm\"uller, M. Hinze, C. Kahle, T. Kiel: {\em A goal-oriented 
dual-weighted adaptive finite element approach for the optimal 
control of a nonsmooth Cahn--Hilliard--Navier--Stokes system.}
WIAS Preprint No. 2311, Berlin 2016, \pier{pp.~1-27.}

\bibitem{HW3}
M. Hinterm\"uller, T. Kiel, D. Wegner: {\em Optimal control of a semidiscrete 
Cahn--Hilliard--Navier--Stokes system with non-matched fluid densities}. 
\pier{SIAM J. Control Optim. {\bf 55} (2017), 1954-1989.}

\bibitem{hw}
M. Hinterm\"uller, D. Wegner: {\em
Distributed optimal control of the Cahn--Hilliard system 
including the case of a double-obstacle homogeneous free energy density}. 
SIAM J. Control Optim. {\bf 50} (2012), 388-418. 

\bibitem{HW1}
M. Hinterm\"uller, D. Wegner: {\em Optimal control of a semidiscrete 
Cahn--Hilliard--Navier--Stokes system}.
SIAM J. Control Optim. {\bf 52} (2014), 747-772.

\bibitem{HW2}  
M.  Hinterm\"uller,  D.  Wegner: {\em Distributed  and  boundary  control  problems  for  the  semidiscrete Cahn--Hilliard/Navier--Stokes system with nonsmooth Ginzburg--Landau energies}. 
Isaac Newton Institute Preprint Series No. NI14042-FRB (2014), \pier{pp.~1-29.}

\bibitem{Kudla}
C. Kudla, A.T. Blumenau, F. B\"ullesfeld, N. Dropka, C. Frank-Rotsch, 
F. Kiessling, O. Klein, P. Lange, W. Miller,
U. Rehse, U. Sahr, M. Schellhorn, G. Weidemann, M. Ziem, 
G. Bethin, R. Fornari, M. M\"uller, J. Sprekels,
V. Trautmann, P. Rudolph: 
{\em Crystallization of 640 kg mc-silicon ingots under traveling magnetic
field by using a heater-magnet module.} J. Crystal Growth {\bf 365} (2013), 54-58.

\bibitem{Lun}
A. Lunardi:
``Analytic semigroups and optimal regularity in parabolic problems''. 
Birkh\"auser Verlag, Basel (1995).

\bibitem{RS}
E. Rocca, J. Sprekels: {\em Optimal distributed control of a nonlocal 
convective Cahn--Hilliard equation by the velocity in three dimensions}. 
SIAM J. Control Optim. {\bf 53} (2015), 1654-1680.

\bibitem{Simon}
J. Simon: {\em Compact sets in the space $L^p(0,T; B)$}.
Ann. Mat. Pura Appl. \pier{(4)}
{\bf 146} (1987), 65-96.

\bibitem{WaNa} Q.-F. Wang, S.-i. Nakagiri: 
{\em Weak solutions of Cahn--Hilliard equations 
having forcing terms and optimal control problems}. 
Mathematical models in functional equations (Japanese) (Kyoto, 1999), 
S\={u}rikaisekikenky\={u}sho K\={o}ky\={u}roku No. 1128 (2000), 172--180.

\bibitem{ZL1}
X.P. Zhao, C.C. Liu: {\em Optimal control of the convective 
Cahn--Hilliard equation}. Appl. Anal. {\bf 92} (2013), 1028-1045.

\bibitem{ZL2} 
X.P. Zhao, C.C. Liu: {\em Optimal control of the convective 
Cahn--Hilliard equation in 2D case}. Appl. Math. Optim. {\bf 70} (2014), 61-82.

\bibitem{ZW}
J. Zheng, Y. Wang: {\em Optimal control problem for Cahn--Hilliard 
equations with state constraint}.
J. Dyn. Control Syst. {\bf 21} (2015), 257-272. 

\End{thebibliography}

}

\End{document}

%
%

Uniqueness

\gian{We assume now that $\tauO=\tauG$, in addition.
Then, \eqref{deftau}~becomes $\tauO=\tauG=\tau$,
and we notice that the function $(p+\tau q,\pG+\tau\qG)$ is $\calV$-valued.}
Since the problem \Adjoint\ is linear, it is sufficient to prove that the unique solution
with $(\phi_3,\phi_4,\phi_5,\phi_6)=(0,0,0,0)$ is $(p,\pG,q,\qG)=(0,0,0,0)$.
In the next lines, $C_i$,~$i=1,2,\dots$,
denote positive constants that depend only on the structural assumptions 
and the $L^\infty$-norms of $\ub$, $\psi$ and~$\psiG$.

We test \eqref{primaA}, written at the time~$s$, by $(p+\tau q,\pG+\tau\qG)(s)$ 
and integrate over $(t,T)$ with respect to~$s$.
At the same time, we test \eqref{secondaA} by $-(q,\qG)$ and integrate over time as before.
Then, we add the \juerg{resulting} identities.
Since four terms cancel out, by rearranging a little, we \juerg{find} that
\Bsist
  && \frac 12 \, \normaHH{(p+\tau q,\pG+\tau\qG)(t)}^2
  \non
  \\
  && \quad {}
  + \bintQt |q|^2
  + \bintSt |\qG|^2
  + \tau \bintQt |\nabla q|^2
  + \tau \bintSt |\nablaG\qG|^2
  \non
  \\
  && = \bintQt \ub \cdot \nabla p \, (p+\tau q)
  - \bintQt \psi q (p+\tau q)
  - \bintSt \psiG \qG (\pG+\tau\qG)\,.
  \label{pertest1}
\Esist 
\juerg{We} now estimate the \rhs.
The sum of the last two terms is easily handled with the help of the Young inequality as~follows:
\Bsist
  && - \bintQt \psi q (p+\tau q)
  - \bintSt \psiG \qG (\pG+\tau\qG)
  \non
  \\
  && \leq \frac 12 \bintQt |q|^2
  + \frac 12 \bintSt |\qG|^2
  + C_1 \int_t^T \normaHH{(p+\tau q,\pG+\tau\qG)(s)}^2 \, ds \,.
  \non
\Esist
Next, we estimate the transport term by owing to 
\eqref{hpuA} and to the \Holder, Young and Sobolev inequalities.
If we denote by $\CO$ the corresponding Sobolev constant, we have~that
\Bsist
  && \bintQt \ub \cdot \nabla p \, (p+\tau q)
  \leq \int_t^T \norma{\ub(s)}_3 \, \norma{\nabla p(s)}_2 \, \norma{(p+\tau q)(s)}_6 \, ds
  \non
  \\
  && \leq \norma\ub_{\L\infty{\Lx3}} \int_t^T \norma{\nabla p(s)}_2 \, 
    \CO \, \bigl(
      \norma{(p+\tau q)(s)}_2 
      + \norma{\nabla p(s)}_2 
      + \tau \norma{\nabla q(s)}_2
    \bigr)
  \non
  \\
  && \leq \tau \bintQt |\nabla q|^2
  + C_2 \bintQt |\nabla p|^2
  + C_3 \int_t^T \normaHH{(p+\tau q,\pG+\tau\qG)(s)}^2 \, ds \,.
  \non
\Esist
By combining these estimates with \eqref{pertest1}, we obtain
\Bsist
  && \frac 12 \, \normaHH{(p+\tau q,\pG+\tau\qG)(t)}^2
  + \frac 12 \bintQt |q|^2
  + \frac 12 \bintSt |\qG|^2
  + \tau \bintSt |\nablaG\qG|^2
  \non
  \\
  && \leq C_2\bintQt |\nabla p|^2
  + C_4 \int_t^T \normaHH{(p+\tau q,\pG+\tau\qG)(s)}^2 \, ds \,.
  \label{test1}
\Esist
Now, we test \eqref{secondaA} by $(p,\pG)$ and integrate  over~$(t,T)$.
By owing to the Young inequality once more, we have~that
\Bsist
  && \bintQt |\nabla p|^2
  + \bintSt |\nablaG\pG|^2
  \non
  \\
  && = \bintQt q (p + \tau q)
  + \bintSt \qG (\pG + \tau \qG)
  - \tau \bintQt |q|^2
  - \tau \bintSt |\qG|^2
  \non
  \\
  \separa
  && \leq \frac 1{4\tau} \bintQt |p + \tau q|^2
  + \frac 1{4\tau} \bintSt |\pG + \tau \qG|^2
  \non
  \\
  && = \frac 1{4\tau} \int_t^T \normaHH{(p+\tau q,\pG+\tau\qG)(s)}^2 \, ds \,.
  \label{test2}
\Esist
At this point, we add \eqref{test1} and \eqref{test2} multiplied by $C_2$ to each other.
We obtain~that
\Bsist
  && \frac 12 \, \normaHH{(p+\tau q,\pG+\tau\qG)(t)}^2
  + \frac 12 \bintQt |q|^2
  + \frac 12 \bintSt |\qG|^2
  + \tau \bintSt |\nablaG\qG|^2
  + C_2 \bintSt |\nablaG\pG|^2
  \non
  \\
  && \leq C_5 \int_t^T \normaHH{(p+\tau q,\pG+\tau\qG)(s)}^2 \, ds \,.
  \non
\Esist
By applying the Gronwall lemma,
we conclude that $p+\tau q=q=0$ \aeQ.
This implies that $(p,\pG,q,\qG)=(0,0,0,0)$,
and the proof is complete.\QED
 